\documentclass{amsart}

\usepackage{latexsym,amsfonts}
 
\theoremstyle{plain}
\newtheorem{theorem}{Theorem}[section]
\newtheorem{lemma}{Lemma}[section]
\newtheorem{proposition}{Proposition}[section]
 
\theoremstyle{definition}
\newtheorem{definition}[theorem]{Definition}
 
\numberwithin{equation}{section}

\newcommand{\qbin}[2]{\genfrac{[}{]}{0pt}{}{#1}{#2}}
\newcommand{\qbins}[2]{{\textstyle\genfrac{[}{]}{0pt}{}{#1}{#2}}}
\newcommand{\Integer}{\mathbb{Z}}
\newcommand{\II}{\text{\it II}}
\newcommand{\A}{\text A}

\begin{document}

\title{An A$_2$ Bailey Lemma and Rogers--Ramanujan-Type Identities}
 
\author[G.~E.~Andrews]{George~E.~Andrews}
\address{Department of Mathematics, The Pennsylvania State University,
University Park, Pennsylvania 16802, U.S.A.}
\email{andrews@math.psu.edu}
 
\author[A.~Schilling]{Anne Schilling}
\address{Instituut voor Theoretische Fysica, Universiteit van Amsterdam,
Valckenierstraat 65, 1018 XE Amsterdam, The Netherlands}
\email{schillin@phys.uva.nl}
\thanks{The second author is supported by the 
``Stichting Fundamenteel Onderzoek der Materie''.}
\author[S.~O.~Warnaar]{S.~Ole Warnaar}
\address{Instituut voor Theoretische Fysica, Universiteit van Amsterdam,
Valckenierstraat 65, 1018 XE Amsterdam, The Netherlands}
\email{warnaar@phys.uva.nl}
\thanks{The third author is supported by a fellowship of the Royal 
Netherlands Academy of Arts and Sciences.}

\keywords{A$_2$ Bailey lemma, Rogers--Ramanujan identities} 
\subjclass{Primary 05A30, 05A19; Secondary 33D90, 33D15, 11P82}
\date{}

\begin{abstract}
Using new $q$-functions recently introduced by Hatayama {\em et al.}
and by (two of) the authors, we obtain an A$_2$ version of the
classical Bailey lemma. We apply our result, which is distinct from
the A$_2$ Bailey lemma of Milne and Lilly, to derive Rogers--Ramanujan-type
identities for characters of the W$_3$ algebra.
\end{abstract}
 
\maketitle


\section{Introduction}
The celebrated Rogers--Ramanujan identities 
\cite{Rogers94,RR19,Schur17} are given by
\begin{equation}\label{RR1}
1+\sum_{n=1}^{\infty} \frac{q^{n^2}}{(1-q)(1-q^2)\cdots(1-q^n)}=
\prod_{n=1}^{\infty}\frac{1}{(1-q^{5n-1})(1-q^{5n-4})}
\end{equation}
and
\begin{equation}\label{RR2}
1+\sum_{n=1}^{\infty} \frac{q^{n(n+1)}}{(1-q)(1-q^2)\cdots(1-q^n)}=
\prod_{n=1}^{\infty}\frac{1}{(1-q^{5n-2})(1-q^{5n-3})}
\end{equation}
for $|q|<1$.
The fame of these identities lies not only in their beauty and fascinating
history~\cite{Hardy40,Andrews76}, but also in their 
relevance to the theory of partitions and many other branches of 
mathematics and physics. In particular, 
MacMahon~\cite{MacMahon16} and Schur~\cite{Schur17} independently noted that 
the left-hand 
side of \eqref{RR1} is the generating function for partitions into parts 
with difference at least two while the right-hand side generates
partitions into parts congruent to $\pm 1$ modulo $5$.
Similarly, the left-hand side of \eqref{RR2} is the generating function 
for partitions into parts with difference at least two and no parts
equal to 1, while the right-hand side of \eqref{RR2} generates
partitions into parts congruent to $\pm 2$ modulo $5$.

Over the years many generalizations of both the analytic and the
combinatorial statement of the Rogers--Ramanujan identities
have been found, see e.g., 
refs.~\cite{Gordon61,Andrews74,Bressoud79,Bressoud80,Bressoud80b,Andrews84}.
All the cited analytic generalizations are accessible through the classical,
or A$_1$ Bailey lemma and can thus be classified as 
``A$_1$ Rogers--Ramanujan-type identities''.
(We always mean identities of the ``sum=product'' form when referring to 
Rogers--Ramanujan-type identities.)

The proof of Bailey's lemma relies on a $_3\phi_2$ summation known as the 
$q$-Pfaff--Saalsch\"utz formula.
In refs.~\cite{ML92,ML95} Milne and Lilly used a $_6\phi_5$ sum for A$_{n-1}$ 
basic hypergeometric functions to establish a higher-rank version of Bailey's 
lemma. Though Milne and Lilly's result has been applied to yield many beautiful
A$_{n-1}$ generalizations of well-known basic hypergeometric function 
identities~\cite{ML95}, it is remarkable that it has not led to 
A$_{n-1}$ Rogers--Ramanujan-type identities.
In refs.~\cite{Milne92,Milne94} Milne gives U$(n)$ generalizations
of the Rogers--Ramanujan identities related to partitions with
differences between parts at least $n$, but in contrast to the $n=2$ case
Milne's identities are not of Rogers--Ramanujan type.

In this paper a first step is taken towards generalizing the
Rogers--Ramanujan identities to A$_{n-1}$.
After an introduction of the classical Bailey lemma 
(sec.~\ref{secA1} and \ref{secBL})
we prove two summation formulas which give
rise to an A$_2$ Bailey lemma and Bailey chain (sec.~\ref{secA2}).
One of these summation formulas involves new $q$-functions recently introduced 
in refs.~\cite{HKKOTY98,Kirillov98,SW98}.
In the second part of the paper applications of the
A$_2$ Bailey lemma are presented. First we prove several Rogers--Ramanujan-type
identities for characters of the W$_3$ algebra, which generalize the
Rogers--Ramanujan identities \eqref{RR1} and \eqref{RR2} and their extensions
due to Gordon and Andrews (sec.~\ref{secsi}--\ref{secRR}).
As a second application we give an A$_2$ generalization of Bressoud's
identities for partitions with even moduli (sec.~\ref{secM33k}).
Finally our lemma is applied to yield summation formulas
for Kostka polynomials (sec.~\ref{secKostka}).

\subsection{Notation}
Throughout the paper the following notations will be used.
The $q$-shifted factorial is defined for all integers $n$ by
\begin{equation*}
(a;q)_{\infty}=(a)_\infty=\prod_{k=0}^{\infty}(1-aq^k)
\end{equation*}
and
\begin{equation*}
(a;q)_n=(a)_n=\frac{(a)_{\infty}}{(aq^n)_{\infty}}.
\end{equation*}
Note in particular that $1/(q)_n=0$ when $n<0$.
We also use the condensed notation of Gaspar and Rahman~\cite{GR90}, i.e.,
\begin{equation*}
(a_1,a_2,\dots,a_m;q)_n=(a_1,a_2,\dots,a_m)_n=(a_1)_n(a_2)_n\dots(a_m)_n.
\end{equation*}
The Gaussian polynomial or $q$-binomial coefficient is defined as
\begin{equation*}
\qbin{n}{m}=\begin{cases}
\displaystyle
\frac{(q)_n}{(q)_m(q)_{n-m}} & \text{for } 0\leq m \leq n \\[2mm]
0 & \text{otherwise.}\end{cases}
\end{equation*}
At times we find it convenient to display both lower entries of the
$q$-binomial, writing $\qbins{m+n}{m,n}$.
We often use the basic hypergeometric function notation
\begin{equation*}
_{r+1}\phi_r \Bigl[\genfrac{}{}{0pt}{}{a_1,a_2,\dots,a_{r+1}}
{b_1,\dots,b_r};q,z\Bigr]=
\sum_{n=0}^{\infty}
\frac{(a_1,a_2,\dots,a_{r+1})_n}
{(q,b_1,\dots,b_r)_n}z^n.
\end{equation*}
Finally we set $\binom{n}{2}=n(n-1)/2$ for $n\in \Integer$, and	
adopt the convention that $\sum_n$ stands for a sum over all 
integers $n$.

\section{The A$_1$ Bailey lemma}\label{secA1}
One of the most elegant approaches towards proving $q$-series identities
is provided by the Bailey lemma. In trying to catch the essence of Rogers' 
second proof of the Rogers--Ramanujan identities~\cite{Rogers17}, 
Bailey~\cite{Bailey49} was led to the following definition.
\begin{definition}\label{defBP}
A pair of sequences $\alpha=\{\alpha_L\}_{L\geq 0}$ and
$\beta=\{\beta_L\}_{L\geq 0}$ that satisfies 
\begin{equation}\label{BP}
\beta_L = \sum_{r=0}^L \frac{\alpha_r}{(q)_{L-r}(aq)_{L+r}}
\end{equation}
forms a Bailey pair relative to $a$.
\end{definition}
With this definition we have the following lemma, known as Bailey's
lemma, but in its present form due to Andrews~\cite{Andrews84} 
and (in special cases) to Paule~\cite{Paule85}.
\begin{lemma}\label{BaileyLemma}
Let $(\alpha,\beta)$ form a Bailey pair relative to $a$. Then the pair
$(\alpha',\beta')$ given by
\begin{equation}\label{BL1}
\begin{aligned}
\alpha'_L &=
\frac{(\rho_1)_L (\rho_2)_L (aq/\rho_1 \rho_2)^L}
{(aq/\rho_1)_L (aq/\rho_2)_L} \: \alpha_L \\[2mm]
\beta'_L &= \sum_{r=0}^L
\frac{(\rho_1)_r (\rho_2)_r (aq/\rho_1 \rho_2)^r
(aq/\rho_1 \rho_2)_{L-r}}
{(aq/\rho_1)_L (aq/\rho_2)_L (q)_{L-r}} \: \beta_r
\end{aligned}
\end{equation}
again forms a Bailey pair relative to $a$.
\end{lemma}
The strength of this lemma is that it can be iterated leading to
the Bailey chain 
$$(\alpha,\beta)\to (\alpha',\beta')\to (\alpha'',\beta'')\to \cdots .$$
An important special case of Bailey's lemma, and the one that will be 
generalized in section~\ref{secA2}, follows when $\rho_1,\rho_2$ tend to 
infinity. Using
\begin{equation*}
\lim_{a\to \infty} a^{-n} (a)_n = (-1)^n q^{\binom{n}{2}},
\end{equation*}
equation~\eqref{BL1} reduces to
\begin{equation}\label{BL2}
\alpha'_L = a^L q^{L^2} \alpha_L, \qquad
\beta'_L = \sum_{r=0}^L \frac{a^r q^{r^2}}{(q)_{L-r}}\:\beta_r.
\end{equation}
The simplest and perhaps most important application follows by
taking the Bailey pair~\cite{Andrews84}
\begin{equation}\label{43}
\alpha_L =\frac{(1-aq^{2L})(a)_L(-1)^L q^{\binom{L}{2}}}{(1-a)(q)_L},
\qquad \beta_L =\delta_{L,0}
\end{equation}
which is an immediate consequence of a simple $_4 \phi_3$ summation 
(ref.~\cite{GR90}, (2.3.4)).
Iterating this Bailey pair using \eqref{BL2} and setting $a=q^{\ell}$,
$\ell=0,1$, yields the identities
\begin{equation}\label{AGfinite}
\sum_r \frac{(-1)^r q^{((2k+1)r-1)r/2}a^{rk}}{(q)_{L-r}(q)_{L+r+\ell}}=
\sum_{n_1,\dots,n_{k-1}}
\frac{a^{n_1+\cdots+n_{k-1}} q^{n_1^2+\cdots+n_{k-1}^2}}
{(q)_{L-n_1}(q)_{n_1-n_2} \ldots (q)_{n_{k-2}-n_{k-1}}(q)_{n_{k-1}}}.
\end{equation}
Letting $L$ tend to infinity and using Jacobi's triple product identity
(ref.~\cite{GR90}, (II.28)) gives ($|q|<1$)
\begin{equation}\label{AG1}
\sum_{n_1,\dots,n_{k-1}}
\frac{q^{n_1^2+\cdots+n_{k-1}^2}}
{(q)_{n_1-n_2} \ldots (q)_{n_{k-2}-n_{k-1}}(q)_{n_{k-1}}}
=\frac{(q^k,q^{k+1},q^{2k+1};q^{2k+1})_{\infty}}
{(q)_{\infty}}
\end{equation}
and
\begin{equation}\label{AG2}
\sum_{n_1,\dots,n_{k-1}}
\frac{q^{n_1^2+\cdots+n_{k-1}^2+n_1+\cdots+n_{k-1}}}
{(q)_{n_1-n_2} \ldots (q)_{n_{k-2}-n_{k-1}}(q)_{n_{k-1}}}
=\frac{(q,q^{2k},q^{2k+1};q^{2k+1})_{\infty}}
{(q)_{\infty}}.
\end{equation}
For $k=2$ these are the Rogers--Ramanujan identities \eqref{RR1} and 
\eqref{RR2}, whereas for $k\geq 3$ they are identities of 
Andrews~\cite{Andrews74}, related to Gordon's partition 
theorem~\cite{Gordon61}.

In the remainder of this paper we frequently use Bailey pairs such as those
of equation \eqref{43}, but unlike \eqref{43} they are generally very 
cumbersome to write down explicitly.
We therefore adopt the practice of not writing down Bailey pairs explicitly,
but to only present (polynomial) identities which imply Bailey pairs.
The latter are generally much more compact.
When $a=q^{\ell}$, the Bailey pair \eqref{43} can, for example, be caught 
in the identity
\begin{equation}\label{alt}
\sum_r (-1)^r q^{r(r-1)/2}\qbin{2L+\ell}{L-r}
=\frac{q^{L(L+\ell)}(q)_{\ell}}{(q^{-\ell})_{-L}(q)_{-L}},
\qquad \ell=0,1,\dots 
\end{equation}
for $2L+\ell\geq 0$. When $L\geq 0$ the right-hand side simplifies to
$(q)_{\ell}\delta_{L,0}$ and one can (after some work) recover \eqref{43}.
For general $2L+\ell\geq 0$ equation \eqref{alt} follows by taking $x=1$ in 
the $q$-binomial formula (ref.~\cite{GR90}, (II.4))
\begin{equation*}
\sum_r (-1)^r x^r q^{r(r-1)/2}\qbin{2L+\ell}{L-r}
=(x)_L(q/x)_{L+\ell}.
\end{equation*}
In section~\ref{secsi} we will consider an A$_2$ generalization of \eqref{alt}
that implies A$_2$ versions of \eqref{AG1} and \eqref{AG2}. 

\section{Bypassing the Bailey lattice}\label{secBL}
Identities \eqref{AG1} and \eqref{AG2} are the $i=k$ and $i=1$ instances
of~\cite{Andrews74}
\begin{equation}\label{AG3}
\sum_{n_1,\dots,n_{k-1}}
\frac{q^{n_1^2+\cdots+n_{k-1}^2+n_i+\cdots+n_{k-1}}}
{(q)_{n_1-n_2} \ldots (q)_{n_{k-2}-n_{k-1}}(q)_{n_{k-1}}}
=\frac{(q^i,q^{2k-i+1},q^{2k+1};q^{2k+1})_{\infty}}
{(q)_{\infty}}
\end{equation}
true for all $i=1,\dots,k$ and $|q|<1$.
In ref.~\cite{AAB87}, Agarwal, Andrews and Bressoud derive the above result
using an extension of the Bailey chain known as the Bailey lattice 
(see also ref.~\cite{Bressoud88}).
In section~\ref{secap} we wish to derive the A$_2$ counterpart of the
identities \eqref{AG3}, but we are faced with the problem that we do 
have an A$_2$ Bailey chain but not a lattice.

Here we describe a little trick which permits the derivation of \eqref{AG3} 
for all $i$ from just the $a=1$ Bailey pair of equation~\eqref{43} 
and the Bailey chain, bypassing the need for a Bailey lattice. Later we
apply the same procedure in the A$_2$ setting.
As a first step we again take the Bailey pair of equation \eqref{43} with
$a=1$ and iterate once. This leads to the identity \eqref{AGfinite} 
for $a=1$ and $k=1$. Now the sum is split according to the parity of $r$
giving the following bounded form of Euler's identity:
\begin{equation}\label{euler}
\sum_r \Bigl\{
q^{(6r+1)r}\qbin{2L}{L-2r}-q^{(2r+1)(3r+1)}\qbin{2L}{L-2r-1}\Bigr\}
=\frac{(q)_{2L}}{(q)_L}.
\end{equation}
This can be rewritten as
\begin{equation}\label{euler2}
\sum_r \frac{q^{(6r+1)r}(1-q^{4r+1})}{(q)_{L-2r}(q)_{L+2r+1}}
=\frac{1}{(q)_L}
\end{equation}
which is an equation that implies a Bailey pair relative to $a=q$.
Iterating $k-i$ times leads to
\begin{equation*}
\sum_r q^{2r(2r+1)(k-i)}
\frac{q^{(6r+1)r}(1-q^{4r+1})}{(q)_{L-2r}(q)_{L+2r+1}}
=\sum_{n_1,\dots,n_{k-i}}
\frac{q^{n_1^2+\cdots+n_{k-i}^2+n_1+\cdots+n_{k-i}}}
{(q)_{L-n_1}\cdots (q)_{n_{k-i-1}-n_{k-i}}(q)_{n_{k-i}}}.
\end{equation*}
Now rewrite this again in the form of equation \eqref{euler}, and combine the 
two terms in the summand of the left-hand side into one to arrive back at a 
form similar to \eqref{AGfinite}. Explicitly,
\begin{equation*}
\sum_r \frac{(-1)^r q^{((2k-2i+3)(r+1)-2)r/2}}{(q)_{L-r}(q)_{L+r}}
=\sum_{n_1,\dots,n_{k-i}}
\frac{q^{n_1^2+\cdots+n_{k-i}^2+n_1+\cdots+n_{k-i}}}
{(q)_{L-n_1} \ldots (q)_{n_{k-i-1}-n_{k-i}}(q)_{n_{k-i}}}
\end{equation*}
and we are back to an equation which implies a Bailey pair relative to $a=1$.
Iterating $i-1$ times finally gives
\begin{equation*}
\sum_r \frac{(-1)^r q^{((2k+1)(r+1)-2i)r/2}}{(q)_{L-r}(q)_{L+r}}
=\sum_{n_1,\dots,n_{k-1}}
\frac{q^{n_1^2+\cdots+n_{k-1}^2+n_i+\cdots+n_{k-1}}}
{(q)_{L-n_1} \ldots (q)_{n_{k-2}-n_{k-1}}(q)_{n_{k-1}}}
\end{equation*}
so that, letting $L$ tend to infinity and using the Jacobi triple product 
identity, we obtain \eqref{AG3}. 

At the heart of the above sequence of steps was the rewriting of \eqref{euler}
into \eqref{euler2} and, later, the reverse of this. This rewriting comes down 
to
\begin{equation*}
\frac{q^{-2r}}{(q)_{L-2r}(q)_{L+2r}}-
\frac{q^{2r+1}}{(q)_{L-2r-1}(q)_{L+2r+1}}=
\frac{q^{-2r}-q^{2r+1}}{(q)_{L-2r}(q)_{L+2r+1}}.
\end{equation*}
As we wish to generalize to higher rank it is essential
to recognize the above equality as the $B_1=-2r+1$, $B_2=2r+2$ case
of the determinant evaluation
\begin{equation*}
\det_{1\leq i,j\leq 2}\biggl(\frac{q^{j(j-B_i)}}{(q)_{L+B_i-j}}\biggr)=
\frac{q^{5-B_1-2B_2}(1-q^{B_2-B_1})}{(q)_{L+B_1-1}(q)_{L+B_2-1}}.
\end{equation*}

\section{An A$_2$ Bailey lemma}\label{secA2}
\subsection{Motivation}
In order to generalize Bailey's lemma we return to 
definition~\ref{defBP} of the classical, or A$_1$ Bailey pair.
Setting $a=q^{\ell}$ with $\ell=0,1,\dots$, we can rewrite \eqref{BP} as
\begin{equation}\label{bprew}
\beta_L=\frac{1}{(aq)_{2L}} 
\sum_{\substack{k_1 \geq k_2 \\[1mm] k_1+k_2=0}} 
\qbin{2L+\ell}{L-k_1,L-k_2+\ell} \alpha_k
\end{equation}
where $k=(k_1,k_2)$ and where we have chosen to display both lower entries
of the $q$-binomial coefficient.
We now generalize the definition of the A$_1$ Bailey pair by replacing the 
$q$-binomial with its higher-rank analogues.
Recall that the $q$-binomial is the $q$-deformation of the ordinary binomial
coefficient, which can be defined through the expansion
\begin{equation}\label{bin}
(x_1+x_2)^L = \sum_{\lambda_1,\lambda_2} x_1^{\lambda_1} x_2^{\lambda_2}
\binom{L}{\lambda_1,\lambda_2}.
\end{equation} 
The $r$th elementary symmetric function in $n$ variables 
is defined as~\cite{Macdonald95}
\begin{equation*}
e_r(x_1,\dots,x_n)=
\sum_{1\leq i_1<i_2<\dots<i_r\leq n} x_{i_1} x_{i_2} \dots x_{i_r}.
\end{equation*}
Identifying the left-hand side of \eqref{bin} as $(e_1(x_1,x_2))^L$ we 
consider the following A$_{n-1}$ generalization of the binomial coefficient
\begin{equation*}
\prod_{a=1}^{n-1} \bigl(e_a(x_1,\dots,x_n)\bigr)^{L_a}=
\sum_{\lambda} x_1^{\lambda_1}\dots x_n^{\lambda_n}
\binom{L_1,\dots,L_{n-1}}{\lambda_1,\dots,\lambda_n}.
\end{equation*}
Since $e_r$ is homogeneous of degree $r$ one sees that
$\binom{L_1,\dots,L_{n-1}}{\lambda_1,\dots,\lambda_n}=0$ if 
$\sum_{i=1}^n \lambda_i
\neq \sum_{j=1}^{n-1}j L_j$.
It is straightforward to write down an explicit expression for the generalized 
binomial as
\begin{equation}\label{genrank}
\binom{L_1,\dots,L_{n-1}}{\lambda_1,\dots,\lambda_n}
=\sum_r \frac{\prod_{a=1}^{n-1} L_a!}
{\prod_{a=1}^{n-1} \prod_{1\leq i_1<\dots<i_a\leq n} r_{i_1i_2\dots i_a}},
\end{equation}
where the summation over $r$ denotes a sum over the $r_{i_1\dots i_a}$ 
($a=1,\dots,n-1)$ such that
\begin{align*}
\sum_{a=1}^{n-1} \:
\sum_{\substack{1\leq i_1<\dots<i_a\leq n \\i_b=p\:(1\leq b\leq a)}} 
r_{i_1i_2\dots i_a}=\lambda_p, & \qquad p=1,\dots,n, \\[1mm]
\sum_{1\leq i_1<\dots<i_a\leq n} r_{i_1i_2\dots i_a}=L_a, & 
\qquad a=1,\dots,n-1.
\end{align*}

In refs.~\cite{HKKOTY98} and \cite{SW98} a $q$-deformation of \eqref{genrank}
was introduced, which we call (following the terminology of ref.~\cite{SW98}) 
(completely antisymmetric) A$_{n-1}$ supernomial.
In particular, Hatayama {\em et al.}~\cite{HKKOTY98} 
(see also \cite{Kirillov98}) propose the following 
representation of the A$_{n-1}$ supernomial.
Let $\nu^{(n)}$ denote the conjugate of the partition 
$(1^{L_1}2^{L_2}\dots (n-1)^{L_{n-1}})$, i.e., 
$\nu^{(n)}_j=L_j+\cdots+L_{n-1}$. 
Then, for $\sum_i \lambda_i=|\nu^{(n)}|$,
\begin{equation}\label{supH}
\qbin{L_1,\dots,L_{n-1}}{\lambda_1,\dots,\lambda_n}
=\sum_{\nu}\prod_{a=1}^{n-1}\prod_{j=1}^a
\qbin{\nu_j^{(a+1)}-\nu_{j+1}^{(a+1)}}{\nu_j^{(a)}-\nu_{j+1}^{(a+1)}},
\end{equation}
where the sum over $\nu$ denotes a sum over sequences
$\emptyset=\nu^{(0)}\subset\nu^{(1)}\subset\dots\subset\nu^{(n)}$
of Young diagrams such that each skew diagram $\nu^{(a)}-\nu^{(a-1)}$ is 
a horizontal $\lambda_a$-strip (see for these notions ref.~\cite{Macdonald95}).

If we now restrict \eqref{supH} to $n=3$, and set
$\nu^{(1)}=\lambda_1$, $\nu^{(2)}=(\lambda_1+m,\lambda_2-m)$ and 
$\nu^{(3)}=(L_1+L_2,L_2)$ we get
\begin{equation}\label{defsupern1}
\qbin{L_1,L_2}{\lambda_1,\lambda_2,\lambda_3}=\sum_m
\qbin{\lambda_1-\lambda_2+2m}{m}\qbin{L_1}{\lambda_1-L_2+m}
\qbin{L_2}{\lambda_2-m}
\end{equation}
for $\lambda_1+\lambda_2+\lambda_3=L_1+2L_2$ and zero otherwise.
In subsequent manipulations we also employ 
the following representation of the A$_2$ supernomial, which, 
unlike \eqref{defsupern1}, is a manifest $q$-deformation of the $n=3$ case 
of \eqref{genrank},
\begin{equation}\label{defsupern2}
\qbin{L_1,L_2}{\lambda_1,\lambda_2,\lambda_3}=\sum_r 
\frac{q^{r_1 r_{23}}(q)_{L_1}(q)_{L_2}}
{(q)_{r_1}(q)_{r_2}(q)_{r_3}(q)_{r_{12}}(q)_{r_{13}}(q)_{r_{23}}},
\end{equation}
where the summation over $r$ denotes a sum over $r_1,\dots,r_{23}$ such that
\begin{equation*}
r_1+r_{12}+r_{13}=\lambda_1, \quad
r_2+r_{12}+r_{23}=\lambda_2, \quad
r_3+r_{13}+r_{23}=\lambda_3
\end{equation*}
and
\begin{equation*}
r_1+r_2+r_3=L_1, \quad r_{12}+r_{13}+r_{23}=L_2.
\end{equation*}
To pass from \eqref{defsupern2} to \eqref{defsupern1} one has to 
apply the $q$-Chu--Vandermonde summation (ref.~\cite{GR90}, (II.7))
\begin{equation}\label{qCV1}
_2\phi_1\Bigl[\genfrac{}{}{0pt}{}{q^{-n},b}{c};q,cq^n/b\Bigr]
=\frac{(c/b)_n}{(c)_n}.
\end{equation}

To conclude our discussion of the supernomials we point out a slightly
different viewpoint, see e.g., refs.~\cite{Kirillov95,Kirillov98}.
In this approach the ordinary $q$-binomial is recognized as 
\begin{equation*}
\qbin{L}{\lambda_1,\lambda_2}=\sum_{\eta\,\vdash |\lambda|}
K_{\eta\lambda}K_{\eta'\mu}(q),
\end{equation*}
where $\mu=(1^L)$ and $\lambda=(\lambda_1,\lambda_2)$ a composition such 
that $|\lambda|=L$.
The $K_{\lambda\mu}$ and $K_{\lambda\mu}(q)$ are the Kostka number and Kostka 
polynomial, respectively~\cite{Macdonald95}.
Given this expression it becomes natural to define a generalized
$q$-binomial as
\begin{equation}\label{altdef}
\qbin{L_1,\dots,L_{n-1}}{\lambda_1,\dots,\lambda_n}=
\sum_{\eta\,\vdash |\lambda|}K_{\eta\lambda}K_{\eta'\mu}(q),
\end{equation}
where now $\mu=(1^{L_1}2^{L_2}\dots (n-1)^{L_{n-1}})$ and 
$\lambda=(\lambda_1,\dots,\lambda_n)$ a composition such that 
$|\lambda|=|\mu|$.
The generalized $q$-binomial defined this way coincides with \eqref{supH}.
{}From its definition in \eqref{altdef} it is clear that the (completely
antisymmetric) A$_{n-1}$ supernomials are the connection coefficients
between the elementary symmetric functions and the 
Hall--Littlewood polynomials in $n$ variables~\cite{Macdonald95},
\begin{equation*}
e_{\lambda}(x_1,\dots,x_n)=\sum_{\mu\,\vdash |\lambda|}
\qbin{L_1,\dots,L_{n-1}}{\lambda_1,\dots,\lambda_n}P_{\mu}(x_1,\dots,x_n;q).
\end{equation*}

\subsection{An A$_2$ Bailey lemma}\label{secmain}
We now come to the main result of this paper, a new Bailey lemma for
the algebra A$_2$.
Guided by definition~\ref{defBP} of an A$_1$ Bailey pair and
its rewriting as \eqref{bprew}, we use the supernomial \eqref{defsupern2} 
to define the following A$_2$ Bailey pair.
\begin{definition}[A$_2$ Bailey pair of type I]
Denote $\alpha_k=\alpha_{k_1,k_2,k_3}$ and $\beta_L=\beta_{L_1,L_2}$ and let
$\alpha=\{\alpha_k\}_{\substack{k_1\geq k_2\geq k_3\\[0.5mm] k_1+k_2+k_3=0}}$ 
and $\beta=\{\beta_L\}_{L_1,L_2\geq 0}$ be a pair of sequences that satisfies
\begin{equation*}
\beta_L = \sum_{\substack{k_1\geq k_2\geq k_3 \\[1mm] k_1+k_2+k_3=0}} 
\alpha_k \; \sum_r \frac{q^{r_1 r_{23}}}
{(q)_{r_1}(q)_{r_2}(aq)_{r_3}(q)_{r_{12}}(q)_{r_{13}}(q)_{r_{23}}},
\end{equation*}
where $\sum_r$ denotes a sum over $r_1,\dots,r_{23}$ such that
\begin{equation}\label{r1}
r_1+r_{12}+r_{13}=L_2-k_1, \quad
r_2+r_{12}+r_{23}=L_2-k_2, \quad
r_3+r_{13}+r_{23}=L_2-k_3
\end{equation}
and
\begin{equation}\label{r2}
r_1+r_2+r_3=2L_1-L_2, \quad r_{12}+r_{13}+r_{23}=2L_2-L_1.
\end{equation}
Then $(\alpha,\beta)$ forms an A$_2$ Bailey pair of type I
relative to $a$.
\end{definition} 
Note that for $a=q^{\ell}$, $\ell=0,1,\dots$, the above definition can be 
recast in the form
\begin{equation*}
\beta_L=\frac{1}{(aq)_{2L_1-L_2}(q)_{2L_2-L_1}}
\sum_{\substack{k_1\geq k_2\geq k_3 \\[1mm] k_1+k_2+k_3=0}} 
\qbin{2L_1-L_2+\ell,2L_2-L_1}{L_2-k_1,L_2-k_2,L_2-k_3+\ell} \alpha_k,
\end{equation*}
which is indeed a generalization of \eqref{bprew} using the A$_2$ 
supernomial \eqref{defsupern2}.

As a second definition we need
\begin{definition}[A$_2$ Bailey pair of type II]\negthinspace
Denote $\alpha_k\negthinspace=\alpha_{k_1,k_2,k_3}$ 
and $\beta_L\negthinspace=\beta_{L_1,L_2}$ and let
$\alpha=\{\alpha_k\}_{\substack{k_1\geq k_2\geq k_3\\[0.5mm] k_1+k_2+k_3=0}}$ 
and $\beta=\{\beta_L\}_{L_1,L_2\geq 0}$ be a pair of sequences that satisfies
\begin{equation*}
\beta_L = \sum_{\substack{k_1\geq k_2\geq k_3 \\[1mm] k_1+k_2+k_3=0}} 
\alpha_k \;
\frac{(aq)_{L_1+L_2}}
{(aq)_{L_1+k_1}(aq)_{L_1+k_2}(q)_{L_1+k_3}
(q)_{L_2-k_1}(q)_{L_2-k_2}(aq)_{L_2-k_3}}.
\end{equation*}
Then $(\alpha,\beta)$ forms an A$_2$ Bailey pair of type II
relative to $a$.
\end{definition} 
With these definitions our A$_2$ version of Bailey's lemma reads
\begin{theorem}\label{theoremBLA2}
Let $(\alpha,\beta)$ form an A$_2$ Bailey pair of type T=I,II relative to $a$.
Then the pair $(\alpha',\beta')$ given by
\begin{equation}\label{BLA2}
\begin{aligned}
\alpha'_k &= a^{k_1+k_2} q^{\frac{1}{2}(k_1^2+k_2^2+k_3^2)} \alpha_k \\[2mm]
\beta'_L &= f_L^{(T)} \sum_{r_1=0}^{L_1}\sum_{r_2=0}^{L_2}
\frac{a^{r_1} q^{r_1^2-r_1 r_2+r_2^2}}
{(q)_{L_1-r_1}(q)_{L_2-r_2}} \: \beta_r 
\end{aligned}
\end{equation}
forms an A$_2$ Bailey pair of type II relative to $a$.
Here $f_L^{(I)}=(aq)^{-1}_{L_1+L_2}$ and $f_L^{(\II)}=1$.
\end{theorem}
Iteration of theorem~\ref{theoremBLA2} gives the A$_2$ Bailey chains
\begin{equation*}
(\alpha,\beta)_I \to (\alpha',\beta')_{\II} 
\to (\alpha'',\beta'')_{\II} \to \cdots
\end{equation*}
and
\begin{equation*}
(\alpha,\beta)_{\II} \to (\alpha',\beta')_{\II} 
\to (\alpha'',\beta'')_{\II} \to \cdots.
\end{equation*}

\begin{proof}[Proof of theorem~\ref{theoremBLA2}]
The transformation of a Bailey pair of type I to a Bailey pair of type II
follows from the summation formula
\begin{multline}\label{sum1}
\sum_{L_1=0}^{M_1} 
\sum_{L_2=0}^{M_2} 
\frac{a^{L_1} q^{L_1^2-L_1 L_2+L_2^2}}
{(q)_{M_1-L_1}(q)_{M_2-L_2}} 
\sum_r \frac{q^{r_1 r_{23}}}
{(q)_{r_1}(q)_{r_2}(aq)_{r_3}(q)_{r_{12}}(q)_{r_{13}}(q)_{r_{23}}} \\
= \frac{a^{k_1+k_2}q^{\frac{1}{2}(k_1^2+k_2^2+k_3^2)}(aq)^2_{M_1+M_2}}
{(aq)_{M_1+k_1}(aq)_{M_1+k_2}(q)_{M_1+k_3}
(q)_{M_2-k_1}(q)_{M_2-k_2}(aq)_{M_2-k_3}},
\end{multline}
where $k_1+k_2+k_3=0$ and with the same restrictions 
on the sum over $r$ as given by equations \eqref{r1} and \eqref{r2}.
We begin by rewriting the sum over $r$ as a sum over just $r_{12}$ and
$r_{13}$ using these restrictions and by then making the shifts
$r_{12}\to r_{12}+L_2$ and $r_{13}\to r_{13}+L_2$.
After replacing $r_{12}$ by $m$ and $r_{13}$ by $n$ 
the left-hand side of \eqref{sum1} thus becomes
\begin{align*}
\sum_{m,n} q^{(k_1+m+n)(m+n)} &
\sum_{L_1=0}^{M_1}
\frac{a^{L_1} q^{L_1(L_1+k_1+m+n)}}
{(q)_{M_1-L_1}(q)_{L_1-k_2+n}(aq)_{L_1+m-k_3}(q)_{-L_1-m-n}} \notag \\
\times &
\sum_{L_2=0}^{M_2}
\frac{q^{L_2(L_2+m+n)}}
{(q)_{M_2-L_2}(q)_{-L_2-k_1-m-n}(q)_{L_2+m}(q)_{L_2+n}}.
\end{align*}
The sums over $L_1$ and $L_2$ have completely decoupled and can both be
carried out using the $q$-Chu--Vandermonde sum \eqref{qCV1}.
Collecting terms we are lucky again as the remaining double sum has once 
more decoupled, leaving us to sum
\begin{align*}
(aq)_{M_1+k_1}(q)_{M_2-k_1}  &
\sum_m\frac{q^{m(m+k_1+k_2)}}
{(aq)_{M_1+m-k_3}(q)_{M_2+m}(q)_{-k_2-m}(q)_{-k_1-m}} \notag \\
\times & \sum_n\frac{a^{k_2-n}q^{n(n+k_1+k_3)}}
{(q)_{M_1+n-k_2}(q)_{M_2+n}(aq)_{-k_3-n}(q)_{-k_1-n}}.
\end{align*}
Now the other instance of $q$-Chu--Vandermonde is needed
(ref.~\cite{GR90}, (II.6))
\begin{equation}\label{qCV2}
_2\phi_1\Bigl[\genfrac{}{}{0pt}{}{q^{-n},b}{c};q,q\Bigr]
=\frac{b^n (c/b)_n}{(c)_n}
\end{equation}
allowing both summations to be performed. After a few simplifications 
this yields the right-hand side of \eqref{sum1}.

The transformation of a type II to a type II Bailey pair follows from 
\begin{multline}\label{sum2}
\sum_{L_1=0}^{M_1}\sum_{L_2=0}^{M_2} 
\frac{a^{L_1}q^{L_1^2-L_1 L_2+L_2^2}}
{(q)_{M_1-L_1}(q)_{M_2-L_2}} \\
\times \frac{(aq)_{L_1+L_2}}
{(aq)_{L_1+k_1}(aq)_{L_1+k_2}(q)_{L_1+k_3}
(q)_{L_2-k_1}(q)_{L_2-k_2}(aq)_{L_2-k_3}} \\[1mm]
= \frac{a^{k_1+k_2}q^{\frac{1}{2}(k_1^2+k_2^2+k_3^2)}(aq)_{M_1+M_2}}
{(aq)_{M_1+k_1}(aq)_{M_1+k_2}(q)_{M_1+k_3}
(q)_{M_2-k_1}(q)_{M_2-k_2}(aq)_{M_2-k_3}}
\end{multline}
for $k_1+k_2+k_3=0$.
Without loss of generality we can assume $k_1\geq k_2 \geq k_3$ so that
$k_1\geq 0$ and $k_3\leq 0$.
Shifting $L_2\to L_2+k_1$ the left-hand side becomes
\begin{multline*}
\sum_{L_1}\frac{a^{L_1}q^{L_1^2-L_1 k_1+k_1^2}}
{(q)_{M_1-L_1}(aq)_{L_1+k_1}(aq)_{L_1+k_2}(q)_{L_1+k_3}} \\
\times \sum_{L_2}\frac{q^{L_2(L_2-L_1+2k_1)}(aq)_{L_1+L_2+k_1}}
{(q)_{M_2-L_2-k_1}(q)_{L_2}(q)_{L_2+k_1-k_2}(aq)_{L_2+k_1-k_3}}.
\end{multline*}
Next we require the transformation
\begin{equation*}
\lim_{c\to\infty} { _3\phi_2 }
\Bigl[\genfrac{}{}{0pt}{}{q^{-n},b,c}{d,e};q,\frac{deq^n}{bc}\Bigr]=
\frac{1}{(e)_n} \; { _2\phi_1 }
\Bigl[\genfrac{}{}{0pt}{}{q^{-n},d/b}{d};q,eq^n\Bigr]
\end{equation*}
which is a special case of the $q$-Kummer--Thomae--Whipple transformation
for $_3\phi_2$ series (ref.~\cite{GR90}, (III.9), see also \eqref{eq9}).
Applying this to the sum over $L_2$ and making some simplifications leads to
\begin{equation*}
\frac{q^{k_1^2}}{(q)_{M_2-k_2}}
\sum_j \frac{q^{j(j+k_1-k_2-k_3)}}{(q)_j(q)_{M_2-k_1-j}(aq)_{k_1-k_3+j}}
\sum_{L_1}\frac{a^{L_1}q^{L_1(L_1-k_1-j)}}
{(q)_{M_1-L_1}(aq)_{L_1+k_2}(q)_{L_1+k_3-j}}.
\end{equation*}
Now shift $L_1\to L_1-k_3+j$ and perform the sum over $L_1$ using
\eqref{qCV1} with $b\to\infty$. This gives
\begin{equation*}
\frac{a^{-k_3}q^{k_1^2-k_2k_3}}{(q)_{M_2-k_2}(aq)_{M_1+k_2}}
\sum_j \frac{a^j q^{j(j+k_1-k_3)}}
{(q)_j(q)_{M_1+k_3-j}(q)_{M_2-k_1-j}(aq)_{k_1-k_3+j}}.
\end{equation*}
The final sum can again be carried out using \eqref{qCV1} yielding the
desired right-hand side of \eqref{sum2}.
\end{proof}

\section{Applications}\label{secap}
In this section several applications of our A$_2$ Bailey lemma will be given.
Many of the resulting Rogers--Ramanujan-type identities can be recognized
as identities for characters of the W$_3$ algebra, and to facilitate
discussions we give a brief introduction to the characters of the
W$_n$ algebras in section~\ref{sec51}. In section~\ref{secsi} a supernomial
identity is proven which serves as a seed to the A$_2$ Bailey lemma. Then,
in sections~\ref{sec33kp1}-\ref{secM33k} A$_2$ Rogers--Ramanujan-type 
identities are derived. Finally, some identities for Kostka
polynomials are proven in section~\ref{secKostka}.

\subsection{Characters of the W$_n$ algebra}\label{sec51}
The W$_n$ algebra, introduced by Zamolodchikov~\cite{Zamolodchikov85}
and Fateev and Lykyanov~\cite{FL88}, is an A$_{n-1}$ 
generalization of the well-known Virasoro algebra.
The minimal series of W$_n$ algebras is labelled by two integers $p$ and $p'$
such that $\gcd(p,p')=1$ and $n\leq p<p'$, with central charge 
\begin{equation*}
c=(n-1)\Bigl(1-\frac{n(n+1)(p'-p)^2}{pp'}\Bigr).
\end{equation*}
In the remainder we denote this series by $M(p,p')_n$. 
Let $\Lambda_0,\dots,\Lambda_{n-1}$ be the fundamental weights
and $P_+(n,\ell)$ the set of level-$\ell$ dominant integral weights of
A$_{n-1}^{(1)}$. 
Then the $M(p,p')_n$ character of the highest-weight representation labelled
by the pair $\xi\in P_+(n,p-n)$ and $a\in P_+(n,p'-n)$ is given 
by~\cite{Mizoguchi91}
\begin{equation*}
\chi_{\xi,a}^{(p,p')}(q)=\frac{1}{\eta(q)^{n-1}}
\sum_{\alpha\in Q}\sum_{w\in {\mathcal W}_n} \epsilon(w)
q^{\frac{1}{2}pp'|\alpha-(\xi+\rho)/p+w(a+\rho)/p'|^2}.
\end{equation*}
Here $Q$ is the root lattice and ${\mathcal W}_n$ the Weyl group of A$_{n-1}$, 
$\epsilon(w)=(-1)^{\ell(w)}$ with $\ell(w)$ the length of $w$, $\rho$ is the
Weyl vector and $\eta(q)$ is the Dedekind $\eta$-function.
We remark that the characters of $M(p,p')_n$ can be identified
with the branching coefficients of the coset pair  
$(\A_{n-1}^{(1)}\oplus \A_{n-1}^{(1)},\A^{(1)}_{n-1})$ at levels
$p/(p'-p)-n$, $1$ and $p'/(p'-p)-n$, respectively~\cite{JMO88,Nakanishi90}.

Of particular interest to us here are the characters of $M(n,k)_n$.
In this case $\xi=0$ and the characters are labelled by a single 
level-$(k-n)$ dominant integral weight which we write as
$j=\sum_{a=0}^{n-1}(j_a-1)\Lambda_a$, with $j_0,\dots,j_{n-1}\geq 1$
and $j_0+\dots+j_{n-1}=k$.
Invoking the A$_{n-1}$ Macdonald identity~\cite{Dyson62,Macdonald72}
\begin{multline}\label{MDI}
\sum_{k_1+\cdots+k_n=0}
\prod_{i=1}^n x_i^{nk_i} q^{\frac{n}{2}k_i^2+ik_i}
\prod_{1\leq i<j\leq n}(1-x_ix_j^{-1}q^{k_i-k_j}) \\
=(q)^{n-1}_{\infty}\prod_{1\leq i<j\leq n}(x_ix_j^{-1},q x_jx_i^{-1})_{\infty}
\end{multline}
the sum over $Q$ and ${\mathcal W}_n$ (i.e., the affine Weyl group) can be 
carried out to yield the following product representation of the (normalized)
$M(n,k)_n$ characters
\begin{equation}\label{Wn}
\chi_j^{(n,k)}(q)=\chi_{(j_0,\dots,j_{n-1})}^{(n,k)}(q)=
\frac{(q^k;q^k)_{\infty}^{n-1}}
{(q;q)_{\infty}^{n-1}}
\prod_{a=1}^{n-1}\prod_{b=0}^{n-1}
(q^{j_b+j_{b+1}+\cdots+j_{a+b-1}};q^k)_{\infty}
\end{equation}
with the convention that $j_i=j_{i'}$ if $i\equiv i' \pmod{n}$.

\subsection{A supernomial identity}\label{secsi}
As a first application of the A$_2$ Bailey lemma we wish to repeat the
working of section~\ref{secBL} to obtain A$_2$ analogues of 
equation~\eqref{AG3}. Our starting point is the following supernomial 
analogue of identity~\eqref{alt}, which is the first instance of 
a hierachy of supernomial identities conjectured in ref.~\cite{SW98}.
\begin{proposition}\label{prop_id}
For $L_1,L_2,\ell$ integers such that
$2L_1-L_2+\ell\geq 0$, $2L_2-L_1\geq 0$ and $\ell\geq 0$,
\begin{multline}\label{superid}
\sum_{k_1+k_2+k_3=0}\sum_{\sigma\in S_3}\epsilon(\sigma)
q^{\frac{1}{2}\sum_{s=1}^3(3k_s-2\sigma_s)k_s}
\qbins{2L_1-L_2+\ell,2L_2-L_1}
{L_2-3k_1+\sigma_1-1+\ell,L_2-3k_2+\sigma_2-2,L_2-3k_3+\sigma_3-3}  \\
=\frac{q^{L_1^2-L_1 L_2+L_2^2+\ell L_1} (q)_{\ell}}
{(q^{-\ell})_{-L_1}(q)_{-L_2}},
\end{multline}
where $S_3$ is the permutation group on $1,2,3$ and
$\epsilon(\sigma)$ is the sign of the permutation $\sigma$.
\end{proposition}

Note that for $L_1,L_2\geq 0$ the right-hand side simplifies to
$(q)_{\ell}\delta_{L_1,0}\delta_{L_2,0}$. Proposition~\ref{prop_id}
is a corollary of the following stronger result.
\begin{lemma}\label{lem_strong}
For $L,a,\ell$ integers such that $0\leq 2a\leq 3L+\ell$ and $\ell\geq 0$,
\begin{multline}\label{strong}
\sideset{}{'}\sum_{\substack{k_1+k_2+k_3=0\\\sigma\in S_3}}
\epsilon(\sigma)
q^{\frac{1}{2}\sum_{s=1}^3(3k_s-2\sigma_s)k_s}
\qbins{3L-2a+\ell,a}{L-3k_1+\sigma_1-1+\ell,L-3k_2+\sigma_2-2,
L-3k_3+\sigma_3-3}\\
=(-1)^t q^{\binom{t+1}{2}+(L-a)(L-a+t)}\qbins{3L-a+\ell}{2L-a+t},
\end{multline}
where the sum on the left-hand side is such that $t=3k_1+1-\sigma_1$ is fixed.
\end{lemma}
Replacing $L\to L_2$ and $a\to 2L_2-L_1$, summing over $t$ using 
the $q$-binomial theorem
$_1\phi_0(q^{-n};-;q,z)=(zq^{-n})_n$ (ref.~\cite{GR90}, (II.4))
and respecting the ranges of $L_1,L_2$ and $\ell$ yields \eqref{superid}.
Hence we are left to prove \eqref{strong}.

For the proof of lemma~\ref{lem_strong} we need the following
three functions introduced by Andrews {\em et al.}~\cite{ABBBFV87},
\begin{equation}\label{gde}
\begin{aligned}
\gamma(N,M)&=\sum_i \Bigl\{ 
q^{6i^2+2i}\qbin{N+M}{N-3i}-q^{6i^2-2i}\qbin{N+M}{N-3i+1}\Bigr\}\\
\delta(N,M)&=\sum_i \Bigl\{ 
q^{6i^2-i}\qbin{N+M}{N-3i}-q^{6i^2-5i+1}\qbin{N+M}{N-3i+1}\Bigr\}\\
\varepsilon(N,M)&=\sum_i \Bigl\{ 
q^{6i^2-4i}\qbin{N+M}{N-3i}-q^{6i^2-8i+2}\qbin{N+M}{N-3i+1}\Bigr\}.
\end{aligned}
\end{equation}
Observe that $\gamma$ and $\varepsilon$ obey the symmetries
\begin{equation}\label{sym}
\begin{aligned}
\gamma(N,M)&=-\gamma(M-1,N+1)\\
\varepsilon(N,M)&=-\varepsilon(M+2,N-2).
\end{aligned}
\end{equation}
All three functions have another sum representation stated in the following 
lemma.
\begin{lemma}\label{lem_gde}
The functions $\gamma,\delta$ and $\varepsilon$ as defined in \eqref{gde}
satisfy
\begin{equation}\label{gde_sum}
\begin{aligned}
\gamma(N,M)&=\sum_{h=0}^{\infty} (-1)^{h+1} q^{(2h+1)(N+1)+\binom{h+1}{2}}
\frac{(q^{M-N-h-1})_{2h+1}}{(q)_{2h+1}}\\
\delta(N,M)&=\sum_{h=0}^{\infty} (-1)^h q^{2hN+\binom{h+1}{2}}
\frac{(q^{M-N-h+1})_{2h}}{(q)_{2h}}\\
\varepsilon(N,M)&=\sum_{h=0}^{\infty} (-1)^h q^{(2h+1)N+\binom{h}{2}-1}
\frac{(q^{M-N-h+2})_{2h+1}}{(q)_{2h+1}}.
\end{aligned}
\end{equation}
\end{lemma}

\begin{proof}
Using the $q$-binomial recurrences
\begin{equation}\label{rec}
\qbin{m+n}{m}=\qbin{m+n-1}{m}+q^n\qbin{m+n-1}{m-1}
=q^m\qbin{m+n-1}{m}+\qbin{m+n-1}{m-1}
\end{equation}
one can show that $\gamma,\delta$ and $\varepsilon$ of equation
\eqref{gde} as well as the expressions in \eqref{gde_sum} obey the
recurrences
\begin{align*} 
\gamma(N,M)&=\gamma(N-1,M)+q^N\delta(N,M-1)\\
\delta(N,M)&=\delta(N-1,M)+q^N\varepsilon(N,M-1)
=\delta(N,M-1)+q^M\gamma(N-1,M)\\
\varepsilon(N,M)&=\varepsilon(N,M-1)+q^M\delta(N-1,M).
\end{align*}
Together with the initial conditions
$\delta(N,N)=\delta(N+1,N)=1$ and $\gamma(N,N+1)=\varepsilon(N,N-2)=0$
this specifies the three functions uniquely. For the representations
in \eqref{gde_sum} these initial conditions can be easily checked.
For $\gamma$ and $\varepsilon$ in \eqref{gde} they follow from the 
symmetries~\eqref{sym}, and for $\delta$ of equation~\eqref{gde} these
are Schur's bounded analogues~\cite{Schur17} of Euler's identity.
\end{proof}

Also necessary for the proof of lemma~\ref{lem_strong} is the next result.

\begin{lemma}\label{lem_id}
For integers $M,b,a,h$ such that
$M\geq 0$, $0\leq b\leq M$, $0\leq a\leq M$ and $0\leq 2h\leq a$,
the following identity holds:
\begin{equation}\label{id}
\sum_{k,m} (-1)^{k+m} q^{\binom{m}{2}-mh+\binom{k+1}{2}+k(b-2h)}
\qbins{M-k-m}{b-m}\qbins{a}{m}\qbins{a-m}{k}
\frac{(q^{h-k-m+1})_{a-2h}}{(q)_{a-2h}}=\qbins{M}{b}\delta_{h,0}.
\end{equation}
\end{lemma}
\begin{proof}
First we establish a symmetry of equation~\eqref{id}
which will be used later. More precisely, we consider the sum on the
left-hand side of the above identity divided by $\qbins{M}{b}$
which we denote, suppressing the $M$ and $h$ dependence, by $g(a,b)$.
(In the following we always assume that $M,b,a,h$ satisfy the ranges
as specified by the lemma.)
 
We start by transforming the sum over $k$ as follows.
When $h-k-m\geq 0$ we apply Sears' $_3\phi_2$ transformation
(ref.~\cite{GR90}, (III.11))
\begin{equation}\label{eq11}
{_3\phi_2}\Bigl[\genfrac{}{}{0pt}{}{q^{-n},b,c}{d,e};q,q\Bigr]=
\frac{(de/bc)_n}{(e)_n} \Bigl(\frac{bc}{d}\Bigr)^n \;
{_3\phi_2}\Bigl[\genfrac{}{}{0pt}{}{q^{-n},d/b,d/c}{d,de/bc};q,q\Bigr],
\end{equation}
with $n=h-m$ and $d=q^{-(M-m)}$.
When $h-k-m<0$ we replace $k\to a-m-k$ and apply Sears' transformation
(ref.~\cite{GR90}, (III.9))
\begin{equation}\label{eq9}
{_3\phi_2}\Bigl[\genfrac{}{}{0pt}{}{q^{-n},b,c}{d,e};q,\frac{deq^n}{bc}\Bigr]=
\frac{(de/bc)_n}{(e)_n} \; {_3\phi_2}
\Bigl[\genfrac{}{}{0pt}{}{q^{-n},d/b,d/c}{d,de/bc};q,eq^n\Bigr],
\end{equation}
with $n=h-1$ and $d=q^{M+m-a-b+1}$, followed by $k\to b-m-k$.
As a result of all this one finds that
\begin{multline*}
g(a,b)=
\sum_{k,m}(-1)^{k+m} q^{\binom{m}{2}-mh+\binom{k+1}{2}+k(a-2h)}
\qbins{M-k-m}{b-k-m}\qbins{a}{m}\qbins{M-a}{k}/\qbins{M}{b} \\
\times \Bigl\{\frac{(q^{b-2h+1})_{h-k-m}}{(q)_{h-k-m}}+
(-1)^b q^{\frac{1}{2}(b-2h)(b-2k-2m+1)}
\frac{(q^{b-2h+1})_{k+m-b+h-1}}{(q)_{k+m-b+h-1}}\Bigr\}.
\end{multline*}
 
Using
\begin{equation*}
(q^{-n})_k=(-1)^k q^{\binom{k}{2}-nk}\frac{(q)_n}{(q)_{n-k}}, \qquad n\geq 0,
\end{equation*}
it can be observed that the two terms in the summand cancel
when $b-2h+1\leq 0$. Hence we can assume $b-2h\geq 0$ and use on
the first term of the summand
\begin{equation*}
\frac{(q^{n+1})_k}{(q)_k}=\frac{(q^{k+1})_n}{(q)_n}, \qquad k,n\geq 0
\end{equation*}
and on the second term
\begin{equation*}
\frac{(q^{n+1})_k}{(q)_k}=(-1)^nq^{\binom{n+1}{2}+nk}
\frac{(q^{-k-n})_n}{(q)_n},\qquad k,n\geq 0.
\end{equation*}
This leads to the identity
\begin{multline*}
g(a,b)= \\
\sum_{k,m}(-1)^{k+m} q^{\binom{m}{2}-mh+\binom{k+1}{2}+k(a-2h)}
\qbins{M-k-m}{b-k-m}\qbins{a}{m}\qbins{M-a}{k}
\frac{(q^{h-k-m+1})_{b-2h}}{(q)_{b-2h}}/\qbins{M}{b}.
\end{multline*}
(The two contributing terms to this equation correspond to $k+m<h+1$ and
$m+k\geq h+1$, respectively).
Reshuffling the four $q$-binomials this can be recognized as $g(b,a)$.
Note that in establishing this symmetry of $g$ we have also
proven \eqref{id} for $b-2h<0$.

We now come to the actual proof of lemma~\ref{lem_id}.
First assume that $h=0$. The third $q$-binomial implies that the summand 
vanishes unless $a-m-k\geq 0$. But in this range 
$(q^{1-k-m})_a$ is nonzero for $m=k=0$ only and the result is immediate.

The case when $h$ is strictly positive is not so straight-forward
and will be proven by induction on $M$ and $b$.
Let us denote the left-hand side of \eqref{id} by $f(M,b,a,h)$.
By~\eqref{rec} $f$ satisfies the recurrence
\begin{equation*}
f(M,b,a,h)=f(M-1,b,a,h)+q^{M-b}f(M-1,b-1,a,h),
\end{equation*}
for $0\leq b\leq M$ and $0\leq a<M$, so that \eqref{id} will be proven
if we can show its validity in the three cases $b=0$, $b=M$ and $a=M$.
By the symmetry between $a$ and $b$ it is sufficient to treat the
cases $b=0$ and $b=M$ only.

The case $b=0$ is the simplest. All we need to show is that
\begin{equation*}
\sum_k (-1)^k q^{\binom{k+1}{2}-2kh}\qbin{a}{k}
\frac{(q^{h-k+1})_{a-2h}}{(q)_{a-2h}}=0,
\end{equation*}
for $0<2h\leq a$. To achieve this we split the summand according to whether
$h-k\geq 0$ or $h-k+1\leq 0$. In the latter case we shift $k\to k-h+a+1$.
Now both terms admit the sum over $k$ to be performed using
the $q$-Chu--Vandermonde summation \eqref{qCV2}, yielding
$(-1)^h q^{-\frac{h}{2}(3h-1)}\qbins{2h-1}{h}$
and its negative, and we are done.

When $b=M$ we are to show that
\begin{equation*}
\sum_m (-1)^m q^{\binom{m}{2}-mh}\qbin{a}{m}
\frac{(q^{h-m+1})_{a-2h}}{(q)_{a-2h}}=0,
\end{equation*}
for $0<2h\leq a$. Again we treat $h-m\geq 0$ and $h-m+1\leq 0$ separately
and shift $m\to m-h+a+1$ in the latter case.
In hypergeometric notation we then need to show the identity
\begin{equation*}
(-1)^h q^{\binom{h+1}{2}}\qbin{a-h}{h} 
{_2\phi_1}\left[\genfrac{}{}{0pt}{}{q^{-h},q^{-a}}{q^{-a+h}};q,q^h\right]
=\qbin{a}{h-1}
{_2\phi_1}\left[\genfrac{}{}{0pt}{}{q^{-h+1},q^{a-2h+1}}
{q^{a-h+2}};q,q^h\right].
\end{equation*}
To make further progress we apply Jackson's transformation 
(ref.~\cite{GR90}, (III.7)).
\begin{equation}\label{Jackson}
{_2\phi_1}\Bigl[\genfrac{}{}{0pt}{}{q^{-n},b}{c};q,z\Bigr]
=\frac{(c/b)_n}{(c)_n}\;
{_3\phi_2}\Bigl[\genfrac{}{}{0pt}{}{q^{-n},b,bzq^{-n}/c}
{bq^{1-n}/c,0};q,q\Bigr].
\end{equation}
Dropping common factors this yields 
\begin{equation*}
{_3\phi_2}\left[\genfrac{}{}{0pt}{}{q^{-h},q^{-h},q^{-a}}{q^{-2h+1},0};
q,q\right]=q^{h(h-a-1)}
{_3\phi_2}\left[\genfrac{}{}{0pt}{}{q^{-h},q^{-h+1},q^{a-2h+1}}
{q^{-2h+1},0};q,q\right],
\end{equation*}
which can be recognized as a specialization of Sears' $_3\phi_2$ 
transformation \eqref{eq11}.
\end{proof}

We finally come to the proof of lemma~\ref{lem_strong}.

\begin{proof}[Proof of lemma \ref{lem_strong}]
Since in equation~\eqref{strong} $\sigma_1$ only takes the values $1,2$ or 
$3$, the condition $t=3k_1+1-\sigma_1$ fixes both $k_1$ and $\sigma_1$. 
Therefore the restrictions $k_1+k_2+k_3=0$ and $\sigma_1+\sigma_2+\sigma_3=6$ 
reduce the primed sum in \eqref{strong} to a sum over only two independent
summation variables. These we choose to be $\tau=\sigma_2-2$ and $i=k_2$.
Further setting $k_1=\eta$
and using the representation \eqref{defsupern1} of the supernomial
with $m\to L-3i-m+\tau$, the left-hand side of \eqref{strong} reads
\begin{multline}\label{t1}
\sum_{i,m,\tau}
\epsilon(\tau) q^{3(\eta^2+i^2+i\eta)+(4-\tau-2\sigma_1)\eta
+(2-2\tau-\sigma_1)i} \\
\times \qbins{3L-2a+\ell}{2L-a+\ell-3i-m-t+\tau}
\qbins{2L+\ell-3i-2m-t+\tau}{L-3i-m+\tau}\qbins{a}{m}.
\end{multline}
Here $\epsilon(\tau)$ denotes the sign of the permutation
$(\sigma_1,\tau+2,4-\sigma_1-\tau)\in S_3$. Now apply
\begin{equation*}
\qbin{A}{B}=\sum_{k=0}^n (-1)^k q^{\binom{k+1}{2}+Bk}\qbin{n}{k}
\qbin{A+n-k}{B+n},
\end{equation*}
(which is a specialization of the $q$-Chu--Vandermonde summation \eqref{qCV2})
with $n=a-m$ to the first $q$-binomial in \eqref{t1}, to get
\begin{multline}\label{t2}
\sum_{m,k} (-1)^k q^{\binom{k+1}{2}+(L-a+\ell-t)k
+(L-m)(L+m-a+t)+\eta(t-\sigma_1+3)}
\qbins{3L-a+\ell-m-k}{L+\ell-m-t}\qbins{a-m}{k}\qbins{a}{m}\\
\times \Bigl\{ \sum_{i,\tau} \epsilon(\tau) 
q^{6i^2+i(2t+6m-3a-1-4\tau)-\tau(2m-a+t-\eta-\tau)}
\qbins{2L-a-k+t}{L-3i-m+\tau}\Bigr\}_{1/q}.
\end{multline}
The expression within curly brackets can be expressed in terms of the functions 
in \eqref{gde} depending on the parity of $a$.
When $a=2\alpha$ one finds
\begin{align*}
q&^{-\frac{1}{2}\{3(\alpha-\eta-m)^2+(2\sigma_1-1)(\alpha-\eta-m)
+(\sigma_1-1)(\sigma_1-2)\}}\\
\times \bigl\{&\chi(\alpha+m+t \text{ even})\;
\delta(L+\frac{1}{2}(t-3\alpha+m),L+\frac{1}{2}(t-\alpha-m)-k)\\
-&\chi(\alpha+m+t \text{ odd})\;
\delta(L+\frac{1}{2}(t-\alpha-m+1)-k,L+\frac{1}{2}(t-3\alpha+m-1))\bigr\}\\
\intertext{and when $a=2\alpha+1$}
q&^{-\frac{1}{2}\{3(\alpha-\eta-m)^2+(2\sigma_1+2)(\alpha-\eta-m)
+(\sigma_1-1)^2\}}\\
\times \bigl\{&\chi(\alpha+m+t \text{ even})\;
\varepsilon(L+\frac{1}{2}(t-3\alpha+m),L+\frac{1}{2}(t-\alpha-m-2)-k)\\
-&\chi(\alpha+m+t \text{ odd})\;q^{-\frac{1}{2}}
\gamma(L+\frac{1}{2}(t-\alpha-m-1)-k,L+\frac{1}{2}(t-3\alpha+m-1))\bigr\}
\end{align*}
where $\chi(\text{true})=1$ and $\chi(\text{false})=0$.

We now insert the representations for $\gamma$,
$\delta$ and $\varepsilon$ provided by lemma~\ref{lem_gde}.
After the shift $h\to \alpha-h$,
expression \eqref{t2} becomes for both parities of $a$
\begin{multline*}
(-1)^t q^{\binom{t+1}{2}+(L-a)(L-a+t)}
\sum_{0\leq h\leq a/2}\: \sum_{k,m}
(-1)^{h+k+m} q^{\binom{h+1}{2}+h(h+2L+t-2a)} \\
\times q^{\binom{m}{2}-mh+\binom{k+1}{2}+k(L+\ell-t-2h)}
\qbins{3L-a+\ell-m-k}{L+\ell-t-m}\qbins{a-m}{k}\qbins{a}{m}
\frac{(q^{h-m-k+1})_{a-2h}}{(q)_{a-2h}}.
\end{multline*}
Applying lemma~\ref{lem_id} to this expression with $M=3L-a+\ell$ and 
$b=L+\ell-t$ (so that indeed $0\leq 2a\leq 3L+\ell$) yields the 
right-hand side of \eqref{strong}.
\end{proof}

\subsection{Character identities for $M(3,3k+1)_3$}\label{sec33kp1}
We now use the supernomial identity~\eqref{superid} to
derive $q$-series identities for $M(3,3k+1)_3$ characters. Since we will 
be dealing with quite lengthy expressions some shorthand notation is needed.
For $v\in \Integer^2$, set $|v|=v_1+v_2$, $(a)_v=(a)_{v_1}(a)_{v_2}$ and 
$v A v=\sum_{i,j=1}^2 v_i A_{i,j} v_j$, with $A$ a two-by-two matrix.
Thus $v_1^2-v_1 v_2+v_2^2=\frac{1}{2}v C v$, with $C$ the Cartan matrix 
of A$_2$. With this notation \eqref{BLA2} can, for example, be rewritten as
\begin{equation*}
\beta'_L = f_L^{(T)} \sum_{r\in \Integer^2_+}
\frac{a^{r_1} q^{\frac{1}{2}r C r}}{(q)_{L-r}} \: \beta_r.
\end{equation*}

We now proceed as in section~\ref{secBL}.
First we take~\eqref{superid} with $\ell=0$ and iterate the corresponding
A$_2$ Bailey pair of type I (relative to $a=1$) once.
This gives a Bailey pair of type II in the form of a doubly bounded
version of the A$_2$ Euler identity,
\begin{equation}\label{A2Euler}
\sum_{k_1+k_2+k_3=0} q^{\frac{3}{2}(k_1^2+k_2^2+k_3^2)}
\sum_{\sigma\in S_3}\epsilon(\sigma)\prod_{s=1}^3
q^{\frac{1}{2}(3k_s-\sigma_s+s)^2-\sigma_s k_s} 
\qbins{L_1+L_2}{L_1+3k_s-\sigma_s+s}=\qbins{L_1+L_2}{L_1}.
\end{equation}
According to our recipe we should rewrite this as a sum over determinants such 
that these determinants can be evaluated explicitly, whilst retaining
an expression of A$_2$ Bailey type.
This can be achieved by making the variable changes $k_i\to -k_{\sigma_i}$
yielding the following determinantal form for the left-hand side of 
\eqref{A2Euler}
\begin{equation}\label{detform}
\sum_{k_1+k_2+k_3=0}q^{6(k_1^2+k_2^2+k_3^2)+4(k_1+2k_2+3k_3)}
\det_{1\leq s,t\leq 3}\Bigl(q^{t(t-s-3k_s)}\qbins{L_1+L_2}{L_1-3k_s+t-s}\Bigr).
\end{equation}
Determinants of the above type have been 
frequently encountered in the theory of plane partitions, and
Krattenthaler has shown that (ref.~\cite{Krattenthaler90}, page~189)
\begin{equation}\label{determinant}
\det_{1\leq s,t\leq n}
\Bigl(q^{t(t-B_s)}\qbins{L_1+L_2}{L_1-B_s+t}\Bigr)
=\prod_{1\leq s<t\leq n}(1-q^{B_t-B_s})
\prod_{s=1}^n \frac{q^{s(s-B_s)}(q)_{L_1+L_2+s-1}}
{(q)_{L_1-B_s+n}(q)_{L_2+B_s-1}}.
\end{equation}
Hence the determinants in \eqref{detform} can be evaluated leading to
\begin{multline*}
\sum_{k_1+k_2+k_3=0}q^{6(k_1^2+k_2^2+k_3^2)+k_1+2k_2+3k_3}
\prod_{1\leq s<t\leq 3}(1-q^{3k_t-3k_s+t-s})
\prod_{s=1}^3\qbins{L_1+L_2+2}{L_1-3k_s-s+3}\\
=\frac{(q)^2_{L_1+L_2+2}}{(q)_{L_1+L_2+1}(q)_{L_1}(q)_{L_2}}.
\end{multline*}
We observe a slight difference with the case treated in section~\ref{secBL}.
There we started with a Bailey pair relative to $a=1$ and rewrote it as
a Bailey pair relative to $a=q$. 
Here the analogous rewriting using determinants has transformed a Bailey pair
relative to $a=1$ into a new Bailey pair relative to 1.
Iterating (the Bailey pair of type II implied by) 
the previous equation $k-i$ times leads to
\begin{multline*}
\sum_{k_1+k_2+k_3=0}
q^{\sum_{s=1}^3 k_s\left(6k_s+s+\frac{3}{2}(k-i)(3k_s+2s)\right)} \\
\times \prod_{1\leq s<t\leq 3}(1-q^{3k_t-3k_s+t-s})
\prod_{s=1}^3\qbins{L_1+L_2+2}{L_1-3k_s-s+3}\\
=\sum_{r^{(j)}\in\Integer^2}
\frac{q^{\sum_{j=1}^{k-i}(
\frac{1}{2}r^{(j)} C r^{(j)}+|r^{(j)}|)}(q)_{L_1+L_2+2}^2}
{(q)_{L-r^{(1)}} \cdots (q)_{r^{(k-i-1)}-r^{(k-i)}}
(q)_{r^{(k-i)}}(q)_{|r^{(k-i)}|+1}}.
\end{multline*}
Now we again apply \eqref{determinant} for $n=3$ so that
\begin{multline*}
\sum_{k_1+k_2+k_3=0}
q^{\frac{1}{2}(3k-3i+4)\sum_{s=1}^3 (3k_s+2s)k_s}
\det_{1\leq s,t\leq 3}\Bigl(
q^{t(t-s-3k_s)}\qbins{L_1+L_2}{L_1-3k_s+t-s}\Bigr) \\
=\sum_{r^{(j)}\in\Integer^2}
\frac{q^{\sum_{j=1}^{k-i} (\frac{1}{2}r^{(j)} C r^{(j)}+|r^{(j)}|)}
(q)_{L_1+L_2}(q)_{L_1+L_2+1}}
{(q)_{L-r^{(1)}} \cdots (q)_{r^{(k-i-1)}-r^{(k-i)}}
(q)_{r^{(k-i)}}(q)_{|r^{(k-i)}|+1}}.
\end{multline*}
As a final step we iterate this equation $i-1$ times arriving at 
\begin{multline}\label{finite5}
\sum_{k_1+k_2+k_3=0}
q^{\frac{1}{2}(3k+1)\sum_{s=1}^3(3k_s+2s)k_s}
\det_{1\leq s,t\leq 3}\Bigl(
q^{it(t-s-3k_s)}\qbins{L_1+L_2}{L_1-3k_s+t-s}\Bigr) \\
=\sum_{r^{(j)}\in\Integer^2}
\frac{q^{\frac{1}{2}\sum_{j=1}^{k-1}
r^{(j)} C r^{(j)}+\sum_{j=i}^{k-1}|r^{(j)}|}
(1-q^{|r^{(i-1)}|+1})(q)^2_{L_1+L_2}}
{(q)_{L-r^{(1)}} \cdots (q)_{r^{(k-2)}-r^{(k-1)}}
(q)_{r^{(k-1)}}(q)_{|r^{(k-1)}|+1}}
\end{multline}
for all $k\geq 1$ and $i=1,\dots,k$. When $i=1$ one should
identify $(1-q^{|r^{(0)}|+1})$ with $(1-q^{L_1+L_2+1})$.

Many more identities can be derived using \eqref{superid}
for general $\ell$. However, for $\ell\geq 1$ one 
cannot evaluate any determinants thus yielding a one-parameter
family of identities for each choice of $\ell(\geq 1)$ only.
A further complication arises from the fact that the right-hand side of
\eqref{superid} is only a simple $(q)_{\ell}\delta_{L_1,0}\delta_{L_2,0}$ for
$\ell=0,1$. As a consequence we find
\begin{multline}\label{genell}
\sum_{k_1+k_2+k_3=0}
q^{\frac{1}{2}(3k+1)\sum_{s=1}^3(3k_s+2s)k_s} \\
\times \det_{1\leq s,t\leq 3}\Bigl(q^{k(t+\ell\chi(t>1))(t-s-3k_s)}
\qbins{L_1+L_2+\ell}{L_1-3k_s+t-s+\ell\chi(t>1)}\Bigr) \\
= \sum_{r^{(j)}\in\Integer^2}
\frac{q^{\sum_{j=1}^k(1+\delta_{j,k})(\frac{1}{2}r^{(j)}C r^{(j)}
+\ell r^{(j)}_1)}
(q)^2_{L_1+L_2+\ell}(q)_{\ell}}
{(q)_{L-r^{(1)}}\cdots (q)_{r^{(k-1)}-r^{(k)}}(q)_{|r^{(k-1)}|+\ell}
(q)_{Cr^{(k)}+\ell e_1}(q^{-\ell})_{-r^{(k)}_1}(q)_{-r^{(k)}_2}},
\end{multline}
where for arbitrary $\ell$ we cannot perform the sum over $r^{(k)}$.
In this equation $\chi(\text{true})=1$, $\chi(\text{false})=0$, 
$Cr^{(k)}=(2r_1^{(k)}-r_2^{(k)},2r_2^{(k)}-r_1^{(k)})$ and $e_1=(1,0)$.
The sum over $r^{(k)}\in\Integer^2$ 
is understood to be a sum such that $Cr^{(k)}+\ell e_1\in \Integer_+^2$.
(Of course $1/(q)_{Cr^{(k)}+\ell e_1}$ is zero outside this range, but we
want to stay clear from $r_1^{(k)}+\ell<0$.)

As mentioned earlier, for $\ell=0,1$ the sum over $r^{(k)}$ trivializes
since the kernel is only non-zero for $r^{(k)}=0$.
This simplifies the right-hand side of \eqref{genell} to
\begin{equation}\label{ell01}
\sum_{r^{(j)}\in\Integer^2}
\frac{q^{\sum_{j=1}^{k-1}(\frac{1}{2}r^{(j)}C r^{(j)}+\ell r^{(j)}_1)}
(q)^2_{L_1+L_2+\ell}}
{(q)_{L-r^{(1)}}(q)_{r^{(1)}-r^{(2)}}\cdots (q)_{r^{(k-1)}}
(q)_{|r^{(k-1)}|+\ell}},
\qquad \ell=0,1.
\end{equation}
For $\ell=2$ we get contributions to the sum for $r^{(k)}=0$ and
$r^{(k)}=-e_1$. Combining these two terms and making the shift
$r^{(j)}\to r^{(j)}-e_1$ $(j=1,\dots,k-1)$ one finds that the right-hand 
side of \eqref{genell} with $\ell=2$ simplifies to
\begin{equation}\label{ell2}
-q^{1-k}\sum_{r^{(j)}\in\Integer^2}
\frac{q^{\sum_{j=1}^{k-1}(\frac{1}{2}r^{(j)}C r^{(j)}+r^{(j)}_2)+r_1^{(k-1)}}
(q)^2_{L_1+L_2+2}}
{(q)_{L+e_1-r^{(1)}}(q)_{r^{(1)}-r^{(2)}}\cdots (q)_{r^{(k-1)}}
(q)_{|r^{(k-1)}|+1}}.
\end{equation}
Note that for $k=2$ equations \eqref{ell01} with $\ell=1$ and \eqref{ell2}
are equal up to an overall factor and a shift in $L_1$.

To transform the above polynomial identities into identities of the
Rogers--Ramanujan type we let $L_1,L_2$ tend to infinity and apply the
Macdonald identity \eqref{MDI} with $n=3$. We can thus claim the following
$q$-series identities.
\begin{theorem}\label{thmAGA2}
Let $|q|<1$ and $k\geq 2$. Then 
\begin{multline}\label{AGA2}
\sum_{r^{(j)}\in\Integer^2}
\frac{q^{\frac{1}{2}\sum_{j=1}^{k-1}
r^{(j)} C r^{(j)}+\sum_{j=i}^{k-1}|r^{(j)}|}
(1-q^{|r^{(i-1)}|+1})}
{(q)_{r^{(1)}-r^{(2)}} \cdots (q)_{r^{(k-2)}-r^{(k-1)}}
(q)_{r^{(k-1)}}(q)_{|r^{(k-1)}|+1}} \\
=\frac{(q^i,q^i,q^{2i},q^{3k-2i+1},q^{3k-i+1},q^{3k-i+1},
q^{3k+1},q^{3k+1};q^{3k+1})_{\infty}}{(q)^3_{\infty}}
\end{multline}
for $i=1,\dots,k$, and
\begin{multline}\label{AGA2b}
\sum_{r^{(j)}\in\Integer^2}
\frac{q^{\sum_{j=1}^{k-1}(\frac{1}{2}r^{(j)}C r^{(j)}+r^{(j)}_1)+
\sigma r_2^{(k-1)}}}
{(q)_{r^{(1)}-r^{(2)}}\cdots (q)_{r^{(k-2)}-r^{(k-1)}}
(q)_{r^{(k-1)}}(q)_{|r^{(k-1)}|+1}} \\
=\frac{(q,q^{k-\sigma},q^{k+1-\sigma},q^{2k+\sigma},q^{2k+1+\sigma},q^{3k},
q^{3k+1},q^{3k+1};q^{3k+1})_{\infty}}{(q)^3_{\infty}}
\end{multline}
for $\sigma=0,1$.
\end{theorem}
Recalling the remark made after \eqref{finite5} we have
$(1-q^{|r^{(0)}|+1})=1$ in the $i=1$ case of \eqref{AGA2}.
Note that despite the factor $(1-q^{|r^{(i-1)}|+1})$ the left-hand
side of~\eqref{AGA2} is a series  with positive coefficients since
the summand is zero for $r^{(k-1)}\not\in\Integer_+^2$ and
$$(1-q^{|r^{(i-1)}|+1})/(q)_{|r^{(k-1)}|+1}=\qbins{|r^{(i-1)}|+1}{1}/
(q^2)_{|r^{(k-1)}|}.$$

Comparing the right-hand sides of the above two formulas with the
character expression \eqref{Wn} of the W$_n$ algebra,
one can identify \eqref{AGA2} and \eqref{AGA2b} as
identities for the $M(3,3k+1)_3$ characters 
$\chi^{(3,3k+1)}_{(i,i,3k+1-2i)}(q)$ and 
$\chi^{(3,3k+1)}_{(1,k-\sigma,2k+\sigma)}(q)$, respectively.
However, in doing so we are confronted with the unpleasant fact
that we have to multiply both sides of 
\eqref{AGA2} and \eqref{AGA2b} by a factor $(q)_{\infty}$.
But then the left-hand sides are no longer series with manifestly 
positive (integer) coefficients.
It thus is desirable to find a summation formula
that allows the left-hand sides of \eqref{AGA2} and \eqref{AGA2b}
to be rewritten in a form that has
an explicit factor $(q)^{-1}_{\infty}$ such that the remaining
expressions are manifestly positive.
Such a summation formula is not known to us generally.

\subsection{A$_2$ analogues of the Rogers--Ramanujan identities}\label{RRA2}
In the following the $k=2$ instance of theorem~\ref{thmAGA2} is
treated in some further detail. Not only because this case warrants special
attention as the A$_2$ generalization of the Rogers--Ramanujan identities, 
but also because one can actually eliminate the spurious $(q)_{\infty}$.
All that is required is
\begin{equation}\label{Jackson2}
\lim_{b\to\infty}{_2\phi_1}\Bigl[\genfrac{}{}{0pt}{}
{q^{-n},b}{c};q,zq^n/b\Bigr]=\frac{1}{(c)_n}\;
{_2\phi_1}\Bigl[\genfrac{}{}{0pt}{}{q^{-n},z/c}{0};q,cq^n\Bigr]
\end{equation}
which is a special case of Jackson's transformation \eqref{Jackson}. 
For example, taking the right-hand side of \eqref{finite5} with $i=k=2$ and 
applying the above transformation to the sum over $r_2$, we find the rewriting
\begin{multline*}
\sum_{r_1,r_2}\frac{q^{r_1^2-r_1 r_2+r_2^2}}
{(q)_{L_1-r_1}(q)_{L_2-r_2}(q)_{r_1}(q)_{r_2}(q)_{r_1+r_2}} \\
=\sum_{r_1,r_2}\frac{q^{r_1^2-r_1 r_2+r_2^2}}
{(q)_{L_1-r_1}(q)_{L_2-r_2}(q)_{L_2+r_1}(q)_{r_1}}\qbin{2r_1}{r_2}.
\end{multline*}
When $L_1,L_2$ tend to infinity the right-hand side has an
extra factor $1/(q)_{\infty}$ as desired.
The two other $M(3,7)_3$ identities (\eqref{finite5} with $i=1$ and $k=2$ and
\eqref{genell} with $\ell=1$ and $k=2$) can be treated similarly.
As a result we have the following A$_2$ versions of the Rogers--Ramanujan
identities
\begin{theorem}[A$_2$ Rogers--Ramanujan identities]
For $|q|<1$,
\begin{equation*}
\sum_{r_1,r_2\geq 0}\frac{q^{r_1^2-r_1 r_2+r_2^2}}{(q)_{r_1}}\qbin{2r_1}{r_2}
=\prod_{n=1}^{\infty}
\frac{1}{(1-q^{7n-1})^2(1-q^{7n-3})(1-q^{7n-4})(1-q^{7n-6})^2},
\end{equation*}
\begin{multline*}
\sum_{r_1,r_2\geq 0}\frac{q^{r_1^2-r_1 r_2+r_2^2+r_1+r_2}}{(q)_{r_1}}
\qbin{2r_1}{r_2} \\
=\prod_{n=1}^{\infty}
\frac{1}{(1-q^{7n-2})(1-q^{7n-3})^2(1-q^{7n-4})^2(1-q^{7n-5})}
\end{multline*}
and
\begin{multline*}
\sum_{r_1,r_2\geq 0}
\frac{q^{r_1^2-r_1 r_2+r_2^2+r_1}}{(q)_{r_1}}\qbin{2r_1+1}{r_2}
=\sum_{r_1,r_2\geq 0}
\frac{q^{r_1^2-r_1 r_2+r_2^2+r_2}}{(q)_{r_1}}\qbin{2r_1}{r_2}\\
=\prod_{n=1}^{\infty}\frac{1}{(1-q^{7n-1})(1-q^{7n-2})(1-q^{7n-3})
(1-q^{7n-4})(1-q^{7n-5})(1-q^{7n-6})}.
\end{multline*}
\end{theorem}
We have not been able to find an identity for the fourth
$M(3,7)_3$ character, corresponding to \eqref{Wn} with $n=3,k=7$ and
$j=(1,3,3)$ (or a permutation thereof).

\subsection{Character identities for $M(3,3k-1)_3$}
Further families of identities follow by making the replacement $q\to 1/q$
in \eqref{detform} and multiplying the resulting equation by
$q^{3L_1 L_2}$. Then proceeding exactly as in 
section~\ref{sec33kp1}, one readily finds that
\begin{multline}\label{finite6}
\sum_{k_1+k_2+k_3=0}
q^{\frac{1}{2}(3k-1)\sum_{s=1}^3(3k_s+2s)k_s}
\det_{1\leq s,t\leq 3}\Bigl(
q^{it(t-s-3k_s)}\qbins{L_1+L_2}{L_1+3k_s+s-t}\Bigr) \\
=\sum_{r^{(j)}\in\Integer^2}
\frac{q^{\frac{1}{2}\sum_{j=1}^{k-1}
r^{(j)} C r^{(j)}+\sum_{j=i}^{k-1}|r^{(j)}|+2r_1^{(k-1)}r_2^{(k-1)}}
(1-q^{|r^{(i-1)}|+1})(q)^2_{L_1+L_2}}
{(q)_{L-r^{(1)}}\cdots (q)_{r^{(k-2)}-r^{(k-1)}}
(q)_{r^{(k-1)}}(q)_{|r^{(k-1)}|+1}}.
\end{multline}
Again, for $i=1$ we must take $(1-q^{|r^{(0)}|+1})=(1-q^{L_1+L_2+1})$.
If we use the $q\to 1/q$ variant of \eqref{genell} with $k=1$
and iterate we find the polynomial identities
\begin{multline}\label{genelldual}
\sum_{k_1+k_2+k_3=0}
q^{\frac{1}{2}(3k-1)\sum_{s=1}^3(3k_s+2s)k_s} \\
\times \det_{1\leq s,t\leq 3}\Bigl(q^{k(t+\ell\chi(t>1))(t-s-3k_s)}
\qbins{L_1+L_2+\ell}{L_1-3k_s+t-s+\ell\chi(t>1)}\Bigr) \\
=\sum_{r^{(j)}\in\Integer^2}
\frac{
(-1)^{\ell}
q^{\sum_{j=1}^k(\frac{1}{2}r^{(j)}C r^{(j)}+\ell r^{(j)}_1)}
(q)_{r^{(k)}_1+\ell}(1)_{r^{(k)}_2}(q)^2_{L_1+L_2+\ell}}
{(q)_{L-r^{(1)}}\cdots (q)_{r^{(k-1)}-r^{(k)}}
(q)_{|r^{(k-1)}|+\ell}
(q)_{Cr^{(k)}+\ell e_1}} \\
\times q^{3r_1^{(k-1)}r_2^{(k-1)}
+\ell(r_2^{(k-1)}-|r^{(k)}|)
-(r_1^{(k-1)}+r_2^{(k)})(r_2^{(k-1)}+r_1^{(k)})-\binom{\ell+1}{2}}.
\end{multline}
For small $\ell$ the sum over $r^{(k)}$ can be performed so that the
right-hand side of \eqref{genelldual} can be replaced by
\begin{equation}\label{ell01dual}
(-q)^{-\ell}\sum_{r^{(j)}\in\Integer^2}
\frac{q^{\sum_{j=1}^{k-1}(\frac{1}{2}r^{(j)}C r^{(j)}+\ell r^{(j)}_1)
+2r_1^{(k-1)}r_2^{(k-1)}+\ell r_2^{(k-1)}}
(q)^2_{L_1+L_2+\ell}}
{(q)_{L-r^{(1)}}(q)_{r^{(1)}-r^{(2)}}\cdots (q)_{r^{(k-1)}}
(q)_{|r^{(k-1)}|+\ell}},
\qquad \ell=0,1,
\end{equation}
and
\begin{equation*}
q^{-2-k}\sum_{r^{(j)}\in\Integer^2}
\frac{q^{\sum_{j=1}^{k-1}(\frac{1}{2}r^{(j)}C r^{(j)}+r^{(j)}_2)
+2r_1^{(k-1)}r_2^{(k-1)}}(q)^2_{L_1+L_2+2}}
{(q)_{L+e_1-r^{(1)}}(q)_{r^{(1)}-r^{(2)}}\cdots (q)_{r^{(k-1)}}
(q)_{|r^{(k-1)}|+1}}, \qquad \ell=2.
\end{equation*}
Letting $L_1,L_2$ tend to infinity and once again using the A$_2$ Macdonald 
identity we obtain further A$_2$ Rogers--Ramanujan-type identities.
\begin{theorem}\label{thmAGA2b}
Let $|q|<1$ and $k\geq 2$. Then
\begin{multline}\label{AGA2c}
\sum_{r^{(j)}\in\Integer^2}
\frac{q^{\frac{1}{2}\sum_{j=1}^{k-1}
r^{(j)} C r^{(j)}+\sum_{j=i}^{k-1}|r^{(j)}|+2r_1^{(k-1)}r_2^{(k-1)}}
(1-q^{|r^{(i-1)}|+1})}
{(q)_{r^{(1)}-r^{(2)}} \cdots (q)_{r^{(k-2)}-r^{(k-1)}}
(q)_{r^{(k-1)}}(q)_{|r^{(k-1)}|+1}} \\
=\frac{(q^i,q^i,q^{2i},q^{3k-2i-1},q^{3k-i-1},q^{3k-i-1},
q^{3k-1},q^{3k-1};q^{3k-1})_{\infty}}{(q)^3_{\infty}}
\end{multline}
for $i=1,\dots,k$, and
\begin{multline}\label{AGA2d}
\sum_{r^{(j)}\in\Integer^2}
\frac{q^{\sum_{j=1}^{k-1}(\frac{1}{2}r^{(j)}C r^{(j)}+r^{(j)}_1)+
2r_1^{(k-1)}r_2^{(k-1)}+\sigma r_2^{(k-1)}}}
{(q)_{r^{(1)}-r^{(2)}}\cdots (q)_{r^{(k-2)}-r^{(k-1)}}
(q)_{r^{(k-1)}}(q)_{|r^{(k-1)}|+1}} \\
=\frac{(q,q^{k-\sigma},q^{k+1-\sigma},q^{2k-2+\sigma},q^{2k-1+\sigma},q^{3k-2},
q^{3k-1},q^{3k-1};q^{3k-1})_{\infty}}{(q)^3_{\infty}}
\end{multline}
for $\sigma=0,1$.
\end{theorem}
In equations \eqref{finite6} and \eqref{ell01dual}--\eqref{AGA2d}
the sum over $r^{(k-1)}$ can be simplified to a one-dimensional sum
by the $q$-Chu--Vandermonde summation \eqref{qCV1} and
transformation \eqref{Jackson2}.
This simplification will be carried out explicitly in the next section
where we deal with the case $k=2$. 
Comparison with \eqref{Wn} shows that
\eqref{AGA2c} and \eqref{AGA2d} are identities for
the $M(3,3k-1)_3$ characters
$\chi^{(3,3k-1)}_{(i,i,3k-2i-1)}(q)$ and
$\chi^{(3,3k-1)}_{(1,k-\sigma,2k-2+\sigma)}(q)$, respectively.

\subsection{The Rogers--Ramanujan identities}\label{secRR}
Again the case $k=2$, corresponding to $M(3,5)_3$, deserves special attention.
First we note that the W$_2$ and W$_3$ algebras $M(2,5)_2$ and $M(3,5)_3$ are
equivalent under level-rank duality~\cite{KNS91}.
Since the $M(2,5)_2$ identities are nothing but the Rogers--Ramanujan
identities \eqref{RR1} and \eqref{RR2}, so should be the $M(3,5)_3$ identities. 
To transform the $k=2$ identities of the previous section
into the Rogers--Ramanujan identities we simplify the two-dimensional
sum over $r^{(1)}$. This can be achieved by observing that for $\sigma=0,1$,
\begin{equation}\label{sim}
\sum_{r_1,r_2}\frac{q^{r_1^2+r_1 r_2+r_2^2+\sigma(r_1+r_2)}}
{(q)_{A-r_1}(q)_{B-r_2}
(q)_{r_1}(q)_{r_2}(q)_{r_1+r_2+\sigma}}
=\sum_r\frac{q^{r(r+\sigma)}}{(q)_{A+B+\sigma}(q)_{A-r}(q)_{B-r}(q)_r}
\end{equation}
by first summing over $r_2$ using the $q$-Chu--Vandermonde summation
\eqref{qCV1} with $b\to\infty$ and then applying transformation
\eqref{Jackson2}.
Now, since $\frac{1}{2}r C r+2r_1 r_2=r_1^2+r_1 r_2+r_2^2$, one can use
\eqref{sim} to rewrite \eqref{finite6} and \eqref{ell01dual} for $k=2$ as
\begin{multline}\label{RRA2rep}
\sum_{k_1+k_2+k_3=0}
 q^{\frac{15}{2}(k_1^2+k_2^2+k_3^2)+5(k_1+2k_2+3k_3)}
\det_{1\leq s,t\leq 3}\Bigl(
q^{(2-\sigma)t(t-s-3k_s)}\qbins{L_1+L_2}{L_1-3k_s+t-s}\Bigr)\\
=\sum_{r\geq 0}\frac{q^{r(r+\sigma)}(q)_{L_1+L_2}}{(q)_{L_1-r}(q)_{L_2-r}(q)_r}
\end{multline}
for $\sigma=0,1$ and
\begin{multline}\label{RRA2repb}
-\sum_{k_1+k_2+k_3=0}
 q^{\frac{15}{2}(k_1^2+k_2^2+k_3^2)+5(k_1+2k_2+3k_3)+1} \\
\times \det_{1\leq s,t\leq 3}\Bigl(
q^{2(t+\chi(t>1))(t-s-3k_s)}\qbins{L_1+L_2+1}{L_1-3k_s+t-s+\chi(t>1)}\Bigr)\\
=\sum_{r\geq 0}\frac{q^{r(r+1)}(q)_{L_1+L_2+1}}{(q)_{L_1-r}(q)_{L_2-r}(q)_r}.
\end{multline}
In the large $L_1,L_2$ limit this reproduces the Rogers--Ramanujan 
identities as it should.
It is intriguing to observe that the level-rank duality, which predicts the 
equality of characters (which are infinite $q$-series), pertains at the 
polynomial level. That is, the doubly bounded versions of the Rogers--Ramanujan
identities \eqref{RRA2rep} and \eqref{RRA2repb} 
admit an A$_1$ representation even at finite $L_1$ and $L_2$. 
Explicitly, one can replace \eqref{RRA2rep} by 
\begin{equation}\label{RRA1rep}
\sum_j (-1)^j q^{\frac{1}{2}j(5j+2\sigma+1)}
\qbins{L_1+L_2}{L_1-j}\qbins{L_1+L_2}{L_2-j}
=\sum_{r\geq 0}\frac{q^{r(r+\sigma)}(q)_{L_1+L_2}}{(q)_{L_1-r}(q)_{L_2-r}(q)_r}
\end{equation}
and \eqref{RRA2repb} by
\begin{equation}\label{RRA1repb}
\sum_j (-1)^j q^{\frac{1}{2}j(5j+3)}
\qbins{L_1+L_2+1}{L_1-j}\qbins{L_1+L_2+1}{L_2-j}
=\sum_{r\geq 0}\frac{q^{r(r+1)}(q)_{L_1+L_2+1}}{(q)_{L_1-r}(q)_{L_2-r}(q)_r}.
\end{equation}
To see that this is indeed correct, take Watson's $_8\phi_7$ transformation 
formula (ref. \cite{GR90}, (III.18))
\begin{multline*}
_8\phi_7 \Bigl[\genfrac{}{}{0pt}{}
{a,qa^{\frac{1}{2}},-qa^{\frac{1}{2}},b,c,d,e,q^{-n}}
{a^{\frac{1}{2}},-a^{\frac{1}{2}},aq/b,aq/c,aq/d,aq/e,aq^{n+1}};q,
\frac{a^2q^{n+2}}{bcde}\Bigr]\\
=\frac{(aq,aq/de)_n}{(aq/d,aq/e)_n} \:
{_4\phi_3} \Bigl[\genfrac{}{}{0pt}{}{q^{-n},d,e,aq/bc}
{aq/b,aq/c,deq^{-n}/a};q,q\Bigr],
\end{multline*}
replace $n\to L_1$ and $e\to q^{-L_2}$ and let $b,c,d\to\infty$. 
Specializing $a$ to $1$ one finds \eqref{RRA1rep} for $\sigma=0$, whereas
specializing $a$ to $q$ yields \eqref{RRA1repb}.
To obtain \eqref{RRA1rep} for $\sigma=1$ we observe that
\begin{multline*}
\sum_{j=-\infty}^{\infty} (-1)^j q^{\frac{1}{2}j(5j+3)}
\qbins{L_1+L_2}{L_1-j}\qbins{L_1+L_2}{L_2-j} \\
=\sum_{j=0}^{\infty} (-1)^j q^{\frac{1}{2}j(5j+3)}
\Bigl(\qbins{L_1+L_2}{L_1-j}\qbins{L_1+L_2}{L_2-j}
-q^{2j+1}\qbins{L_1+L_2}{L_1-j-1}\qbins{L_1+L_2}{L_2-j-1}\Bigr) \\
=\frac{(q)_{L_1+L_2}}{(q)_{L_1+L_2+1}}
\sum_{j=0}^{\infty} (-1)^j q^{\frac{1}{2}j(5j+3)}
(1-q^{2j+1})\qbins{L_1+L_2+1}{L_1-j}\qbins{L_1+L_2+1}{L_2-j}
\end{multline*}
which corresponds to \eqref{RRA1repb}. The equality of the last two lines
follows by application of the determinant identity \eqref{determinant}
with $n=2$, $B_1=1-j$ and $B_2=2+j$.

\subsection{Rogers--Ramanujan-type identities for ``$M(3,3k)_3$''}
\label{secM33k}
As a further application we derive a family of identities generalizing
Bressoud's series~\cite{Bressoud80}
\begin{equation}\label{Bressoudid}
\sum_{n_1,\dots,n_{k-1}}
\frac{q^{n_1^2+\cdots+n_{k-1}^2+n_i+\cdots+n_{k-1}}}
{(q)_{n_1-n_2}\ldots (q)_{n_{k-2}-n_{k-1}}
(q^2;q^2)_{n_{k-1}} } =
\frac{(q^i,q^{2k-i},q^{2k};q^{2k})_{\infty}}{(q)_{\infty}}
\end{equation}
true for all $k\geq 2$, $1\leq i \leq k$ and $|q|<1$.
Comparing the right-hand side with \eqref{Wn} we are led to label 
\eqref{Bressoudid} as identities for the anomalous series $M(2,2k)_2$.
To generalize this to $M(3,3k)_3$, we do not rely on the supernomial identity 
\eqref{superid} as initial condition, but on the following polynomial identity 
of Gessel and Krattenthaler (ref.~\cite{GK97}, (6.18))
\begin{equation*}
\sum_{k_1+k_2+k_3=0} q^{\frac{3}{2}\sum_{s=1}^3 (3k_s+2s)k_s}
\det_{1\leq s,t\leq 3}\Bigl(q^{t(t-s-3k_s)}\qbins{L_1+L_2}{L_1+3k_s+s-t}\Bigr)
=\qbins{L_1+L_2}{L_1}_{q^3}.
\end{equation*}
Clearly this gives rise to a Bailey pair of type II relative to $a=1$.
Along the now clear-trodden path we thus derive for $1\leq i \leq k$,
\begin{multline*}
\sum_{k_1+k_2+k_3=0}
q^{\frac{3}{2}k\sum_{s=1}^3 (3k_s+2s)k_s}
\det_{1\leq s,t\leq 3}\Bigl(
q^{it(t-s-3k_s)}\qbins{L_1+L_2}{L_1+3k_s+s-t}\Bigr) \\
=\sum_{r^{(j)}\in\Integer^2}
\frac{q^{\frac{1}{2}\sum_{j=1}^{k-1}
r^{(j)} C r^{(j)}+\sum_{j=i}^{k-1}|r^{(j)}|}
(1-q^{|r^{(i-1)}|+1})(q^3;q^3)_{|r^{(k-1)}|}(q)_{L_1+L_2}^2}
{(q)_{L-r^{(1)}} \cdots (q)_{r^{(k-2)}-r^{(k-1)}}(q^3;q^3)_{r^{(k-1)}}
(q)_{|r^{(k-1)}|+1}(q)_{|r^{(k-1)}|}}.
\end{multline*}
When $L_1,L_2$ approach infinity this yields the following theorem.
\begin{theorem}
For $|q|<1$, $k\geq 2$ and $i=1,\dots,k$,
\begin{multline*}
\sum_{r^{(j)}\in\Integer^2}
\frac{q^{\frac{1}{2}\sum_{j=1}^{k-1}
r^{(j)} C r^{(j)}+\sum_{j=i}^{k-1}|r^{(j)}|}
(1-q^{|r^{(i-1)}|+1})(q^3;q^3)_{|r^{(k-1)}|}}
{(q)_{r^{(1)}-r^{(2)}}\cdots (q)_{r^{(k-2)}-r^{(k-1)}}(q^3;q^3)_{r^{(k-1)}}
(q)_{|r^{(k-1)}|+1}(q)_{|r^{(k-1)}|}} \\
=\frac{(q^i,q^i,q^{2i},q^{3k-2i},q^{3k-i},q^{3k-i},
q^{3k},q^{3k};q^{3k})_{\infty}}{(q)^3_{\infty}}.
\end{multline*}
\end{theorem}
Comparing with equation \eqref{Wn} this corresponds to identities for
the ``characters'' $\chi^{(3,3k)}_{(i,i,3k-2i)}(q)$.

\subsection{Identities for Kostka polynomials}\label{secKostka}
A different type of application of the A$_2$ Bailey lemma arises from
the representation \eqref{altdef} of the supernomials in terms
of the Kostka polynomials.
Using the summation (ref.~\cite{Fulton97}, page~76)
\begin{equation*}
\sum_{\sigma\in S_n}\epsilon(\sigma)
K_{\eta(\lambda_1+\sigma_1-1,\ldots,\lambda_n+\sigma_n-n)}=\delta_{\eta\lambda}
\end{equation*}
one can invert \eqref{altdef}, resulting in~\cite{SW98}
\begin{equation*}
K_{\lambda'\mu}(q)=\
\sum_{\sigma\in S_n}\epsilon(\sigma)
\qbins{L_1,\dots,L_n}
{\lambda_1+\sigma_1-1,\dots,\lambda_n+\sigma_n-n},
\end{equation*}
where $\mu=(1^{L_1}2^{L_2}\dots (n-1)^{L_{n-1}})$
and $\lambda$ is a partition such that $|\lambda|=|\mu|$.
For $n=3$ this supernomial identity can be used as input to the A$_2$ 
Bailey lemma yielding the following identity for the Kostka polynomials
for $k\geq 1$ and $i=1,\dots,k$,
\begin{multline}\label{itK}
(q)^2_{L_1+L_2}\sum_{r^{(j)}\in\Integer_+^2}
\frac{q^{\frac{1}{2}\sum_{j=1}^k r^{(j)}C r^{(j)}
+\sum_{j=i}^{k-1}|r^{(j)}|}(1-q^{|r^{(i-1)}|+1})K_{\lambda'\mu}(q)}
{(q)_{L-r^{(1)}}\cdots (q)_{r^{(k-1)}-r^{(k)}}
(q)_{|r^{(k-1)}|+1}(q)_{C r^{(k)}}} \\
=q^{\frac{k}{2}(k_1^2+k_2^2+k_3^2)-k(k_1+2k_2+3k_3)}
\det_{1\leq s,t\leq 3}\Bigl(q^{it(t-s+k_s)}\qbins{L_1+L_2}{L_1-k_s+s-t}\Bigr)
\end{multline}
where $\mu=(1^{(Cr^{(k)})_1} 2^{(Cr^{(k)})_2})=(1^{2r^{(k)}_1-r^{(k)}_2}
2^{2r^{(k)}_2-r^{(k)}_1})$
and $\lambda=(|\mu|^3)+(k_1,k_2,k_3)$ with $k_1\geq k_2\geq k_3$, 
$k_1+k_2+k_3=0$.
Also recall that $|r^{(0)}|=L_1+L_2$.
Taking the limit $L_1,L_2\to\infty$ of~\eqref{itK} and evaluating the 
determinant leads to our final theorem.
\begin{theorem}
For $|q|<1$, $k\geq 1$ and $i=1,\dots,k$,
\begin{multline*}
\sum_{r^{(j)}\in\Integer_+^2}
\frac{q^{\frac{1}{2}\sum_{j=1}^k r^{(j)}C r^{(j)}
+\sum_{j=i}^{k-1}|r^{(j)}|}(1-q^{|r^{(i-1)}|+1})K_{\lambda'\mu}(q)}
{(q)_{r^{(1)}-r^{(2)}}\cdots (q)_{r^{(k-1)}-r^{(k)}}
(q)_{|r^{(k-1)}|+1}(q)_{Cr^{(k)}}} \\
=\frac{1}{(q)^3_{\infty}}
q^{\frac{k}{2}(k_1^2+k_2^2+k_3^2)-(k-i)(k_1+2k_2+3k_3)}
\prod_{1\leq s<t\leq 3}(1-q^{i(k_s-k_t-s+t)}).
\end{multline*}
\end{theorem}
The above results can either be viewed as some formal identities for the 
Kostka polynomials, or, by for example taking the representation for the
Kostka polynomials due to Kirillov and Reshetikhin~\cite{KR88} (see
also ref.~\cite{Macdonald95}, page 245), as explicit $q$-series identities.

For $k=1$ the identity \eqref{itK} can be written as
\begin{multline*}
\frac{(q)_{L_1+L_2+2}^2}{(q)_{L_1+L_2+1}}
\sum_{C^{-1}r\in\Integer_+^2}
\frac{q^{\frac{1}{2}rC^{-1}r} K_{\eta\mu}(q)}
{(q)_{L-C^{-1}r}(q)_r} \\
=q^{\frac{1}{2}(\lambda_1^2+\lambda_2^2+\lambda_3^2)}
\prod_{1\leq s<t\leq 3}(1-q^{\lambda_s-\lambda_t-s+t})
\prod_{s=1}^3\qbins{L_1+L_2+2}{L_1-\lambda_s+s-1},
\end{multline*}
where $\mu=(1^{r_1}2^{r_2})$ and 
$\eta=(1^{\lambda_1-\lambda_2}2^{\lambda_2-\lambda_3}3^{\lambda_3+|\mu|/3})$ 
with $\lambda_1\geq \lambda_2\geq \lambda_3$, $|\lambda|=0$.
This is a bounded analogue of the $l=1$, $n=3$ case of 
equation (4.39) of ref.~\cite{HKKOTY98}. Indeed, letting $L_1,L_2$ 
tend to infinity and recognizing 
$$\frac{q^{\frac{1}{2}(\lambda_1^2+\lambda_2^2+\lambda_3^2)}}{(q)^2_{\infty}}
\prod_{1\leq s<t\leq 3}(1-q^{\lambda_s-\lambda_t-s+t})$$
as the branching function $b_{\lambda}(q)$ 
of the level-1 basic representation $V(\Lambda_0)$
of A$_{n-1}^{(1)}$ we find
\begin{equation*}
b_{\lambda}(q)=\sum_{C^{-1}r\in\Integer_+^2}
\frac{q^{\frac{1}{2}rC^{-1}r} K_{\eta\mu}(q)}{(q)_r}
\end{equation*}
of ref.~\cite{HKKOTY98}.


\begin{thebibliography}{99}

\bibitem{AAB87}
A.~K.~Agarwal, G.~E.~Andrews and D.~M.~Bressoud,
{\em The Bailey lattice},
J. Ind. Math. Soc. {\bf 51} (1987), 57--73.

\bibitem{Andrews74}
G.~E.~Andrews,
{\em An analytic generalization of the Rogers--Ramanujan identities
for odd moduli},
Prod. Nat. Acad. Sci. USA {\bf 71} (1974), 4082--4085.

\bibitem{Andrews76}
G.~E.~Andrews,
{\em The Theory of Partitions},
Encyclopedia of Mathematics and its Applications, Vol.~2,
(Addison-Wesley, Reading, Massachusetts, 1976).

\bibitem{Andrews84}
G.~E.~Andrews,
{\em Multiple series Rogers--Ramanujan type identities},
Pacific J. Math. {\bf 114} (1984), 267--283.

\bibitem{ABBBFV87}
G.~E.~Andrews, R.~J.~Baxter, D.~M.~Bressoud, W.~H.~Burge, P.~J.~Forrester
and G.~Viennot,
{\em Partitions with prescribed hook differences},
Europ. J. Combinatorics {\bf 8} (1987), 341--350.

\bibitem{Bailey49}
W.~N.~Bailey,
{\em Identities of the Rogers--Ramanujan type},
Proc. London Math. Soc. (2) {\bf 50} (1949), 1--10.

\bibitem{Bressoud79}
D.~M.~Bressoud,
{\em A generalization of the Rogers--Ramanujan identities for all moduli},
J. Combin. Theory Ser. A {\bf 27} (1979), 64--68.

\bibitem{Bressoud80}
D.~M.~Bressoud,
{\em An analytic generalization of the Rogers--Ramanujan identities
with interpretation},
Quart. J. Maths. Oxford (2) {\bf 31} (1980), 385--399.

\bibitem{Bressoud80b}
D.~M.~Bressoud,
{\em Analytic and combinatorial generalizations of the Rogers--Ramanujan
identities},
Memoirs Amer. Math. Soc. {\bf 24} (1980), 1--54.

\bibitem{Bressoud88}
D.~M.~Bressoud,
{\em The Bailey lattice: An introduction},
in {\em Ramanujan Revisited}, pp. 57--67,
G.~E.~Andrews {\em et al.} eds., (Academic Press, New York, 1988).

\bibitem{Dyson62}
F.~Dyson,
{\em Missed opportunities},
Bull. Amer. Math. Soc. {\bf 78} (1972), 140--156.

\bibitem{FL88}
V.~A.~Fateev and S.~L.~Lykyanov,
{\em The models of two-dimensional conformal quantum field 
theory with $Z_n$ symmetry},
Int. J. Mod. Phys. A {\bf 3} (1988), 507--520.

\bibitem{Fulton97}
W.~Fulton,
{\em Young tableaux: with applications to representation theory and geometry},
London Math. Soc. student texts {\bf 35}, Cambridge University Press (1997).

\bibitem{GR90}
G.~Gasper and M.~Rahman,
{\em Basic Hypergeometric Series},
Encyclopedia of Mathematics and its Applications, Vol.~35, 
(Cambridge University Press, Cambridge, 1990).

\bibitem{GK97}
I.~M.~Gessel and C.~Krattenthaler, 
{\em Cylindric partitions},
Trans. Amer. Math. Soc. {\bf 349} (1997), 429--479.

\bibitem{Gordon61}
B.~Gordon,
{\em A combinatorial generalization of the Rogers--Ramanujan identities},
Amer. J. Math. {\bf 83} (1961), 393--399.

\bibitem{Hardy40}
G.~H.~Hardy, 
{\em Ramanujan},
(Cambridge University Press, London and New York, 1940).

\bibitem{HKKOTY98}
G.~Hatayama, A.~N.~Kirillov, A.~Kuniba, M.~Okado, T.~Takagi and Y.~Yamada,
{\em Character formulae of $\widehat{sl}_n$-modules and inhomogeneous paths},
math.QA/9802085.

\bibitem{JMO88}
M.~Jimbo, T.~Miwa and M.~Okado,
{\em Local state probabilities of solvable lattice models: An A$_{n-1}^{(1)}$
family},
Nucl. Phys. B {\bf 300 [FS22]} (1988), 74--108.

\bibitem{Kirillov95}
A.~N.~Kirillov,
{\em Dilogarithm identities},
Prog. Theor. Phys. Suppl. {\bf 118} (1995), 61--142.

\bibitem{Kirillov98}
A.~N.~Kirillov,
{\em New combinatorial formula for Hall--Littlewood polynomials},
math.QA/\-9803006.

\bibitem{KR88}
A.~N.~Kirillov and N.~Yu.~Reshetikhin,
{\em The Bethe Ansatz and the combinatorics of Young tableaux},
J. Soviet Math. {\bf 41} (1988), 925--955.

\bibitem{Krattenthaler90}
C.~Krattenthaler, 
{\em Generating functions for plane partitions of a given shape}, 
Manuscripta Math. {\bf 69} (1990), 173--202.

\bibitem{KNS91}
A.~Kuniba, T.~Nakanishi and J.~Suzuki,
{\em Ferro-- and antiferro-magnetizations in RSOS models},
Nucl. Phys. B {\bf 356} (1991), 750--774.

\bibitem{Macdonald72}
I.~G.~Macdonald,
{\em Affine root systems and Dedekind's $\eta$-function},
Inv. Math. {\bf 15} (1972), 91--143.

\bibitem{Macdonald95}
I.~G.~Macdonald,
{\em Symmetric Functions and Hall Polynomials}, 2nd edition,
(Clarendon Press, Oxford, 1995).

\bibitem{MacMahon16}
P.~A.~MacMahon,
{\em Combinatory Analysis}, Vol.~2,
(Cambridge University Press, London and New York, 1916).

\bibitem{Milne92}
S.~C.~Milne
{\em Classical partition functions and the $U(n+1)$ Rogers--Selberg identity},
Discrete Math. {\bf 99} (1992), 199-246.

\bibitem{Milne94}
S.~C.~Milne
{\em A $q$-analog of a Whipple's transformation for hypergeometric series 
in $U(n)$},
Adv. in Math. {\bf 108} (1994), 1--76.

\bibitem{ML92}
S.~C.~Milne and G.~M.~Lilly,
{\em The $A_{\ell}$ and $C_{\ell}$ Bailey transform and lemma},
Bull. Amer. Math. Soc. (N.S.) {\bf 26} (1992), 258--263.

\bibitem{ML95}
S.~C.~Milne and G.~M.~Lilly,
{\em Consequences of the $A_{\ell}$ and $C_{\ell}$ Bailey transform and
Bailey lemma},
Discrete Math. {\bf 139} (1995), 319--346.

\bibitem{Mizoguchi91}
S.~Mizoguchi,
{\em The structure of representations for the W$_3$ minimal model},
Int. J. Mod. Phys. A {\bf 6} (1991), 133-162.

\bibitem{Nakanishi90}
T.~Nakanishi,
{\em Non-unitary minimal models and RSOS models},
Nucl. Phys. B {\bf 334} (1990), 745--766.

\bibitem{Paule85}
P.~Paule,
{\em On identities of the Rogers--Ramanujan type},
J. Math. Anal. Appl. {\bf 107} (1985), 225--284.

\bibitem{Rogers94}
L.~J.~Rogers,
{\em Second memoir on the expansion of certain infinite products},
Proc. London Math. Soc. {\bf 25} (1894), 318--343.

\bibitem{Rogers17}
L.~J.~Rogers,
{\em On two theorems of combinatory analysis and some allied identities},
Proc. London Math. Soc. (2) {\bf 16} (1917), 315--336.

\bibitem{RR19}
L.~J.~Rogers and S.~Ramanujan,
{\em Proof of certain identities in combinatory analysis},
Proc. Cambridge Phil. Soc. {\bf 19} (1919), 211--216.

\bibitem{SW98}
A.~Schilling and S.~O.~Warnaar,
{\em Inhomogeneous lattice paths, generalized Kostka polynomials and
A$_{n-1}$ supernomials}, math.QA/9802111.

\bibitem{Schur17}
I.~J.~Schur,
{\em Ein Beitrag zur additiven Zahlentheorie und zur Theorie der
Kettenbr\"uche},
S.-B. Preuss. Akad. Wiss. Phys.-Math. Kl. (1917), 302--321.

\bibitem{Zamolodchikov85}
A.~B.~Zamolodchikov,
{\em Infinite additional symmetries in two-dimensional conformal 
quantum field theory},
Teo. Mat. Fiz. {\bf 65} (1985), 347--359.

\end{thebibliography}
\end{document}